\documentclass[oupthm,crop,info]{CUP-JNL-CJM}
\usepackage{graphicx}
\usepackage{multicol,multirow}
\usepackage{amsmath,amssymb,amsfonts}
\usepackage{ytableau}
\usepackage{tabularray}
\usepackage{mathrsfs}
\usepackage{rotating}
\usepackage{appendix}
\usepackage[numbers]{natbib}
\usepackage{makecell}
\usepackage{subcaption}
\usepackage{mathtools}
\usepackage{stmaryrd}

\usepackage{eucal}
\usepackage{extarrows}
\usepackage{tcolorbox}
\usepackage{tikz-cd}
\usepackage{tkz-graph}
\usetikzlibrary{arrows.meta}
\usepackage{bm}

\theoremstyle{oupplain}
\newtheorem{theorem}{Theorem}[section]

\theoremstyle{oupdefinition}
\newtheorem{definition}{Definition}[section]
\newtheorem{alg}[theorem]{Algorithm}
\theoremstyle{oupremark}
\newtheorem{remark}[theorem]{Remark}
\newtheorem{example}[theorem]{Example}
\theoremstyle{oupproof}

\newtheorem{proofGL}{Proof (general linear group)}
\newtheorem{proofO}{Proof (orthogonal group)}
\newtheorem{proofSp}{Proof (symplectic group)}

\numberwithin{equation}{section}

\articletype{RESEARCH ARTICLE}
\jname{Canadian Mathematical Society}
\jyear{2023}

\newcommand{\C}{\mathbb C}
\newcommand{\N}{\mathbb N}
\DeclareMathOperator{\GL}{GL}
\DeclareMathOperator{\Par}{Par}
\newcommand{\Z}{\mathbb Z}
\DeclareMathOperator{\M}{M}
\DeclareMathOperator{\SM}{SM}
\DeclareMathOperator{\AM}{AM}
\renewcommand{\P}{\mathcal P}
\newcommand{\la}{\lambda}
\DeclareMathOperator{\Hom}{Hom}

\renewcommand{\d}{\partial}
\newcommand{\al}{\alpha}
\newcommand{\be}{\beta}
\renewcommand{\O}{\operatorname{O}}
\DeclareMathOperator{\Sp}{Sp}
\DeclareMathOperator{\SL}{SL}
\renewcommand{\P}{\mathcal P}
\newcommand{\Sd}{{\mathfrak S_d}}
\renewcommand{\c}{c}
\newcommand{\chat}{\widehat{c}}
\renewcommand{\v}[1]{\vcenter{\hbox{#1}}}

\begin{document}

\begin{Frontmatter}

\title[Graphical methods and invariants]{Graphical methods and rings of invariants on the symmetric algebra}

\author{Rebecca Bourn}
\author{William Q. Erickson}
\author{Jeb F. Willenbring}

\authormark{R. Bourn, W. Erickson and J. Willenbring}

\address{\orgname{University of Wisconsin--Milwaukee}, \orgaddress{\city{Milwaukee, WI}, \country{USA}}
\email{bourn@uwm.edu}}
\address{\orgname{Baylor University}, \orgaddress{\city{Waco, TX}, \country{USA}}
\email{will_erickson@baylor.edu}}
\address{\orgname{University of Wisconsin--Milwaukee}, \orgaddress{\city{Milwaukee, WI}, \country{USA}}
\email{jw@uwm.edu}}



\abstract{Let $G$ be a complex classical group, and let $V$ be its defining representation (possibly plus a copy of the dual).
A foundational problem in classical invariant theory is to write down generators and relations for the ring of $G$-invariant polynomial functions on the space $\P^m(V)$ of degree-$m$ homogeneous polynomial functions on $V$.
In this paper, we replace $\P^m(V)$ with the full polynomial algebra $\P(V)$.
As a result, the invariant ring is no longer finitely generated.
Hence instead of seeking generators, we aim to write down linear bases for bigraded components.
Indeed, when $G$ is of sufficiently high rank, we realize these bases as sets of graphs with prescribed number of vertices and edges.
When the rank of $G$ is small, there arise complicated linear dependencies among the graphs,
but we remedy this setback via representation theory:
in particular, we determine the dimension of an arbitrary component in terms of branching multiplicities from the general linear group to the symmetric group.
We thereby obtain an expression for the bigraded Hilbert series of the ring of invariants on $\P(V)$.
We conclude with examples using our graphical notation, several of which recover classical results.}

\end{Frontmatter}

\section{Introduction}
\label{sec:intro}

The motivation for this paper began, in part, with the problem of writing down invariant functions on the space of polynomial-coefficient differential operators (i.e., the Weyl algebra), under the action of the general linear group.  While the goal seems modest, the combinatorial complexity of the answer is certainly not.  The result involves a correspondence between these invariants and unlabeled directed graphs.  In retrospect, it became evident that the linear change of variables for the Weyl algebra had an analogue for the other classical groups.

The archetypal problem in 19th-century invariant theory was to consider the space $S^m(\C^2)$ of binary $m$-forms, where the group $\SL(2, \C)$ acts naturally on $\C^2$, and to describe the ring of polynomial functions on $S^m(\C^2)$  which are invariant under the action of $\SL(2,\C)$.
The problem naturally generalizes to any of the classical groups $\GL_n$, $\O_n$, and $\Sp_{2n}$, where $\C^2$ is replaced by the defining representation $V$ of the group (plus a copy of the dual representation in the case of $\GL_n$).
Since the $m$th symmetric power $S^m(V)$ is finite-dimensional, Hilbert's celebrated basis theorem guarantees that the ring of invariants is finitely generated.
Nearly every paper in early invariant theory contained explicit, often remarkably lengthy computations of generators and relations for rings of invariants.

In this paper, we sum out the dependence on the degree $m$ before looking at the invariants; in other words, we study the invariant functions on the full symmetric (equivalently, polynomial) algebra, rather than on a single graded component.
This idea is related to the classical study of \emph{perpetuants} (or more generally \emph{subvariants}, both terms coined by Sylvester~\cite{SylvesterSubvariants}), which can be viewed as certain $\SL(2,\C)$-semiinvariants of a binary form of infinite order; see also~\cite[p.~326]{GY}, \cite[p.~149]{olver}, \cite{KP}, and references therein.
As a result of our replacing $S^m(V)$ by $S(V)$, the invariant ring is no longer finitely generated.
Hence our goal, unlike the usual goal in classical invariant theory, is not to find generators and relations; rather, we impose a bigradation on the invariant ring, so that each graded component is finite-dimensional, and we aim to write down a linear basis for each component.
In particular, we grade the invariant algebra by degree $d$, in the usual sense for polynomial functions, and by weight $k$ (see Definition~\ref{def:weight}).  
By symmetrizing Hermann Weyl's fundamental theorems of invariant theory, we show that each bigraded component is spanned by the invariants represented by certain graphs (directed, undirected, etc.) with $d$ vertices and $k$ edges.
Our first result (Algorithm~\ref{alg}) gives the explicit map between a graph and its corresponding invariant.
We emphasize that our method easily recovers many classical results, simply by restricting our attention to $m$-regular graphs (i.e., graphs in which every vertex has degree $m$).

Our graphs furnish a true basis whenever the parameters $n,d,k$ fall within a \emph{stable range}, described by a simple inequality relating the rank of the group to $d$ and $k$:
 \begin{equation}
 \label{stable range}
     \GL_n \text{ and } \O_n: \quad n \geq \min\{d,k\}.  \qquad \qquad \Sp_{2n}: \quad n \geq \min\!\left\{\left\lfloor \frac{d}{2}\right\rfloor, \left\lfloor\frac{k}{2}\right\rfloor\right\}.
 \end{equation}
 Hence, for fixed $d$ and $k$, our methods yield a basis for all but finitely many values of $n$.  
 (Also at the  farthest extreme from the stable range, where the defining representation of the group is 1-dimensional, we can easily describe a basis for all $d$ and $k$ in terms of degree sequences of graphs; see Section \ref{sub:n=1}.)

 Outside the stable range, our methods furnish a spanning set rather than a true basis, and therefore we would overshoot the dimension of the graded component by merely counting graphs.  
 The linear dependencies that arise among the graphs are quite complicated, but nonetheless we are able to express the desired dimensions by taking a representation-theoretic approach.
 Our second result (Theorem~\ref{thm:branching}) is a dimension formula for the bigraded components that holds even outside the stable range.  
 The formula is expressed as the sum of branching multiplicities from the general linear group $\GL_d$ to its subgroup $\Sd$ of permutation matrices.
 In this way, we can express the bigraded Hilbert series of the invariant ring.

Although the main problem in this paper seems to be new, the use of graphical methods is nearly as old as invariant theory itself; in fact, it was Sylvester \cite{Sylvester} who coined the term ``graph'' to describe the diagrams he developed in his ``algebro-chemical'' approach to invariants ( see~\cite{Sylvester} and~\cite[p.~366]{GY}).  
Subsequently, Kempe \cite{kempe} and (much later) Olver and Shakiban \cite{OS} improved upon Sylvester's program by associating covariants on binary forms with directed graphs obeying certain syzygies.
Since then, the graphical approach continues to develop; consider, for example, webs and spiders~\cite{Kuperberg} in the invariant theory of low-rank Lie groups.
This development of graphical notation is analogous to similar programs in various fields, most notably the Penrose notation in mathematical physics \cite{Penrose} and Lie theory \cite{Cvitanovic}.   
In fact, the reader may like to begin by skimming the examples in Section~\ref{section:examples}, among which are several classical results translated into our graphical notation.

\section{Preliminaries}

\subsection{Statement of the problem}

Throughout the paper, we let $\M_{nd} \coloneqq \M_{nd}(\C)$ denote the space of complex $n \times d$ matrices, with $\M_n \coloneqq \M_{nn}$.  
By the (complex) \emph{classical groups}, we mean
\begin{itemize} 
\item the \emph{general linear group}
$\GL_n \coloneqq \GL(n,\C) = \{g \in \M_n \mid \det g \neq 0 \}$;
\item 
the \emph{orthogonal group} $\O_n \coloneqq \O(n,\C) = \{ g \in \GL_n \mid g^Tg = I\}$;
\item the \emph{symplectic group}
$\Sp_{2n} \coloneqq \Sp(2n,\C) = \{ g \in \GL_{2n} \mid g^T J g = J\}$,
where $J = \left[\begin{smallmatrix} 0 & I \\ -I & 0 \end{smallmatrix}\right]$.
\end{itemize}
That is, $\O_n$ preserves the symmetric bilinear form $b$ on $\C^n$ given by $b(v,w)=v^Tw$, while $\Sp_{2n}$ preserves the skew-symmetric bilinear form $\omega$ on $\C^{2n}$ given by $\omega(v,w) = v^T J w$.  

In the context of each classical group $G$, we will write $V$ to denote its defining representation, on which it acts naturally by matrix multiplication: 
hence for $\GL_n$ and $\O_n$ we have $V = \C^n$, while for $\Sp_{2n}$ we have $V=\C^{2n}$.  
For $G = \O_n$ or $\Sp_{2n}$, we let $\Psi \coloneqq \P(V)$ denote the space of polynomial functions on $V$; for $G = \GL_n$, we let $\Psi \coloneqq \P(V \oplus V^*)$.  
In each case, $\Psi$ is a representation of $G$ in the natural way, via $g \cdot \psi(v) = \psi(g^{-1}v)$, for $g \in G$, \: $\psi \in \Psi$, and $v \in V$ (or $V \oplus V^*$).  
This, in turn, induces a representation of $G$ on $\P(\Psi)$, the space of polynomial functions on $\Psi$ (explained below), again via $g \cdot f(\psi) = f(g^{-1} \cdot \psi)$ for $f \in \P(\Psi)$.  
Our goal is to describe, as explicitly as possible, the invariant ring
\[
\P(\Psi)^G \coloneqq \Big\{ f \in \P(\Psi) \:\: \Big| \:\: g \cdot f = f \text{ for all }g \in G\Big\}.
\]

\subsection{Notation}
\label{sec:notation}

Let $\N \coloneqq \{0,1,2,\ldots\}$. 
A Greek letter $\alpha = (\al_1,\ldots,\al_n) \in\N^n$ will denote a multi-index  corresponding to the exponent vector of the monomial $x_1^{\al_1} \cdots x_n^{\al_n}$.  
It will be natural to regard $x_i$ as the linear functional on $V$ dual to the $i$th standard basis vector. 
We will use these exponent vectors $\al$ to pick out coefficients of monomials; our desired invariants will then be polynomials in these coefficients.  
Below we organize the details, which vary slightly for each of the three classical groups: 
\begin{itemize}
    \item For $G = \GL_n$:
    \begin{itemize}
        \item We have $\Psi \coloneqq \P(V \oplus V^*) \cong  \C[x_1,\ldots,x_n,\d_1,\ldots,\d_n].$  
        We can regard $\d_i$ as the $i$th standard basis element of the double dual $V^{**}$, although the notation $\d_i$ is indeed meant as an abbreviation for the differential operator $\d/\d x_i$.  
        (See Section \ref{section:examples} for the motivation from the Weyl algebra, i.e., the algebra of polynomial-coefficient differential operators.
        Note, however, that the analysis in this paper will not use the  non-commutative algebra structure of the Weyl algebra, which we view only as a representation of ${\rm GL}_n$.
        )
        \item For $\al,\be \in \N^n$, we write $        \mathbf x^\al \bm{\d}^\be \coloneqq x_1^{\al_1}\cdots x_n^{\al_n} \d_1^{\be_1}\cdots \d_n^{\be_n}$.
        \item We have the monomial basis $\{\mathbf x^\al\bm{\d}^\be \mid \al,\be \in \N^n\}$ for $\Psi$.
        \item We define the linear functional $  \c_{\al,\be} : \Psi \longrightarrow \C$,
    which takes the value 1 on the basis vector $\mathbf x^\al \bm{\d}^\be$ and 0 on all other basis vectors.
    \item Hence we have $\P(\Psi) = \C[\c_{\al,\be}]_{\al,\be \in \N^n}$.
    \end{itemize}
    
\item For $G = \O_n$:

\begin{itemize}
    \item We have $\Psi \coloneqq \P(V) \cong \C[x_1,\ldots,x_n]$.
    \item We have the monomial basis $\{\mathbf x^\al \mid \al \in \N^n\}$ for $\Psi$, where $\mathbf x^\al \coloneqq x_1^{\al_1}\cdots x_n^{\al_n}$.
    \item We define the linear functional $\c_\al: \Psi \longrightarrow \C$, which takes the value 1 on the basis vector $\mathbf x^\al$ and 0 on all other basis vectors.
    \item Hence we have $\P(\Psi) = \C[\c_\al]_{\al\in \N^n}$.
\end{itemize}
\item For $G = \Sp_{2n}$:
\begin{itemize}
    \item We have $\Psi \coloneqq \P(V) \cong \C[x_1,\ldots,x_{2n}]$.
    \item We have the monomial basis $\{\mathbf x^\al \mid \al \in \N^{2n}\}$ for $\Psi$.
    \item We define the linear functional $\c_\al:\Psi \longrightarrow \C$, just as for $\O_n$ above.
    \item Hence we have $\P(\Psi) = \C[\c_\al]_{\al\in \N^{2n}}$.
\end{itemize}
\end{itemize}

The span of the $\c$'s is the graded dual $\Psi^\circ \coloneqq 
\bigoplus_{i=0}^\infty (\Psi_i)^*$, where $\Psi_i$ is the $i$th homogeneous graded component of $\Psi$ with respect to polynomial degree.  
Note that $\Psi^\circ$ is not the same as the full dual space $\Psi^*$, which includes infinite linear combinations of the $\c$'s.  
The ``$\c$'' notation stands for ``coefficient,'' since the $\c$'s extract the coefficient of the corresponding monomial.  
Hence our desired invariants will be polynomials in these $\c$'s, i.e., elements of $\P(\Psi)\cong S(\Psi^\circ)$.  
As is customary in the literature, we define scaled coefficients $\chat$ as follows, writing $\al! \coloneqq \al_1!\cdots \al_n!$:
\[
\chat_{\al,\be} \coloneqq \alpha!\beta! \: \c_{\al,\be} \quad \text{and} \quad \chat_\al \coloneqq \alpha! \: \c_\al.
\]

We will appeal to the graded isomorphism of $G$-modules $\Psi \cong \Psi^\circ$ given by
\begin{equation}    \label{dualbasis}
    \begin{cases}
    \mathbf x^\al\bm{\d}^\be  \longmapsto  \chat_{\al,\be}, & G = \GL_n\\ \phantom{\mathbf x^\al}\mathbf x^\al \longmapsto  \chat_\al, & G = \O_n \text{ or }\Sp_{2n}.
    \end{cases}
\end{equation}
(See \cite[p.~117]{WallachGIT} and \cite[pp.~7--8]{dolgachev}.)  From a different perspective, it is a standard fact (see \cite[Rem.~1.2]{dolgachev}) that the differential operator $\bm{\d}^\alpha$ and the monomial $\mathbf x^\al$ transform dually under a classical group action; hence in each graded component, the isomorphism in \eqref{dualbasis} can equivalently be written as
\begin{equation*}
    \begin{cases}
    \mathbf x^\al\bm{\d}^\be  \longmapsto \bm{\d}^\al \mathbf x^\be, & G = \GL_n\\ \phantom{\mathbf x^\al}\mathbf x^\al \longmapsto  \bm{\d}^\al,& G = \O_n \text{ or }\Sp_{2n}.
    \end{cases}
\end{equation*}
Here we regard $\bm{\d}^\al \mathbf x^\be$ as the operator acting via
$\bm{\d}^\al \mathbf x^\be(\mathbf x^\gamma \bm{\d}^\delta) = \bm{\d}^\al(\mathbf x^\gamma)\cdot \bm{\d}^\delta(\mathbf x^\be)$.  
For example, if $G=\GL_2$ and $\Psi_{4,3} = \operatorname{span}\{\mathbf x^\al \bm{\d}^\be : |\al|=4, \quad |\be|=3\}$, then $\d_1^3 \d_2^{\phantom{1}} x_1^{\phantom{1}} x_2^2$ equals the functional $\chat_{(3,1),(1,2)}$, since it takes the value $3!\cdot 1! \cdot 1! \cdot 2!$ on the monomial $x_1^3 x_2^{\phantom{1}} \d_1^{\phantom{1}}\d_2^2 \in \Psi_{4,3}$, and takes the value 0 on all other monomials in $\Psi_{4,3}$.

In order to obtain finite-dimensional graded components, we impose a bigradation on $\P(\Psi)$.  
First is the natural gradation by \emph{degree}; hence $\P^d(\Psi)$ is the space of homogeneous polynomials of degree $d$ in the $\c$'s.  
The second gradation is by \emph{weight}.  
To make this clear, we write $|\alpha|\coloneqq \sum_i \alpha_i$, and then we define the weight of a monomial in the $\c$'s as follows:
\begin{definition}
\label{def:weight}
    For $G=\GL_n$, the \emph{weight} of the monomial $\c_{\alpha,\beta} \cdots \c_{\psi,\omega}$ equals $\frac{1}{2}(|\alpha| + |\beta| + \cdots + |\psi|+|\omega|)$.
    For $G=\O_n$ or $\Sp_{2n}$, the \emph{weight} of the monomial $\c_\al\! \cdots \c_\omega$ equals $\frac{1}{2}(|\alpha| + \cdots + |\omega|)$.
\end{definition}
As we will see, the factor $\frac{1}{2}$ arises due to the fact that the fundamental invariants are quadratics. 
Now we define $\P_k(\Psi)$ to be the space of polynomials that are homogeneous (i.e., ``isobaric,'' in the classical language) of weight $k$.  
Finally, we set $\P^d_k(\Psi) \coloneqq \P^d(\Psi) \cap \P_k(\Psi)$.

\subsection{Weyl's fundamental theorems}
\label{sub:FFT SFT}

In his monumental book~\cite{weyl}, Hermann Weyl determined generators and relations for the ring of invariants, in the setting where a classical group $G$ acts naturally on an arbitrary number $d$ of vectors (and covectors, when $G = \GL_n$).
For all three classical groups, the first fundamental theorem (FFT) states that the generators are certain quadratic functions $r_{ij}$.
The second fundamental theorem (SFT) states the relations as the determinants (or Pfaffians) of certain minors in the $r_{ij}$.

Combining the FFT and SFT, one obtains an algebra isomorphism $\mu^\sharp$ between the ring of invariants, on one hand, and the coordinate ring of a certain determinantal variety, on the other hand.
(The notation $\mu^\sharp$ is typical in the literature, because this map is the comorphism induced by a matrix multiplication map $\mu$.)
We will use the notation $\SM_d \subset \M_d$ for the subspace of symmetric matrices, and $\AM_d \subset \M_d$ for the subspace of alternating (i.e., skew-symmetric) matrices; we let $z_{ij}$ denote the natural coordinate functions on $\M_d$.  
The aforementioned \emph{determinantal varieties} are denoted by
\begin{align*}
    \M^{\leq n}_d &\coloneqq \{ X \in \M_d \: \mid \operatorname{rank} X \leq n\},\\
    \SM^{\leq n}_d &\coloneqq \{ X \in \SM_d \mid \operatorname{rank} X \leq n\},\\
    \AM^{\leq 2n}_d &\coloneqq \{ X \in \AM_d \mid \operatorname{rank} X \leq 2n\}.
\end{align*}
For our purposes in this paper, we present a version of the combined FFT and SFT below (once for each classical group) in terms of matrix coordinates $x_{ij}$, rather than the coordinate-free presentation used by Weyl.

\begin{theorem}[FFT and SFT for $G = \GL_n$]\
\label{FTs GLn}

Let $G = \GL_n$, acting on $\M_{nd} \oplus \M_{nd} \cong V^{\oplus d} \oplus (V^*)^{\oplus d}$ via $g \cdot (X,Y) = (gX, (g^{-1})^TY)$.
Set
    \begin{equation}
        \label{R GLn}
        R \coloneqq \C[x_{ij}, y_{ij}]_{\substack{1 \leq i \leq n, \\ 1 \leq j \leq d\phantom{,}}} \cong \P(\M_{nd} \oplus \M_{nd}).
    \end{equation}
    Then $R^G$ is generated by the quadratics
    \begin{equation}
        \label{rij GLn}
        r_{ij} \coloneqq \sum_{\ell = 1}^n y_{\ell i} x_{\ell j}, \qquad 1 \leq i,j \leq d.
    \end{equation}
    Moreover, we have the following isomorphism of algebras:
    \begin{align}\label{muGLn}
    \begin{split}
    \mu^{\sharp}:\P(\M_d^{\leq n}) &\longrightarrow R^G,\\
    z_{ij} & \longmapsto r_{ij}.
    \end{split}
    \end{align}
    
\end{theorem}

\begin{theorem}[FFT and SFT for $G = \O_n$]\
\label{FTs On}

Let $G = \O_n$, acting by matrix multiplication on $\M_{nd} \cong V^{\oplus d}$.
Set
    \begin{equation}
        \label{R On}
        R \coloneqq \C[x_{ij}]_{\substack{1 \leq i \leq n, \\ 1 \leq j \leq d\phantom{,}}} \cong \P(\M_{nd}).    \end{equation}
    Then $R^G$ is generated by the quadratics
    \begin{equation}
        \label{rij On}
        r_{ij} \coloneqq \sum_{\ell = 1}^n x_{\ell i} x_{\ell j}, \qquad 1 \leq i \leq j \leq d.
    \end{equation}
    Moreover, we have the following isomorphism of algebras:
    \begin{align}
    \label{muOn}
    \begin{split}
    \mu^{\sharp}:\P(\SM_d^{\leq n}) &\longrightarrow R^G,\\
    z_{ij} & \longmapsto r_{ij}.
    \end{split}
    \end{align}
    
\end{theorem}

\begin{theorem}[FFT and SFT for $G = \Sp_{2n}$]\
\label{FTs Sp2n}

Let $G = \Sp_{2n}$, acting by matrix multiplication on $\M_{2n,d} \cong V^{\oplus d}$.
Set
    \begin{equation}
        \label{R Sp2n}
        R \coloneqq \C[x_{ij}]_{\substack{1 \leq i \leq 2n, \\ 1 \leq j \leq d\phantom{2,}}} \cong \P(\M_{2n,d}).
    \end{equation}
    Then $R^G$ is generated by the quadratics
    \begin{equation}
        \label{rij Sp2n}
        r_{ij} \coloneqq \sum_{\ell = 1}^n (x_{\ell i} x_{\ell+n, j} - x_{\ell + n, i} x_{\ell j}), \qquad 1 \leq i < j \leq d.
    \end{equation}
    Moreover, we have the following isomorphism of algebras:
    \begin{align}\label{muSp2n}
    \begin{split}
    \mu^{\sharp}:\P(\AM_d^{\leq 2n}) &\longrightarrow R^G,\\
    z_{ij} & \longmapsto r_{ij}.
    \end{split}
    \end{align}
    
\end{theorem}

\section{Linear bases in terms of graphs}
\label{section:algorithm}

Weyl's fundamental theorems foreshadow the effectiveness of graphical data in invariant theory.
Indeed, by viewing each $r_{ij}$ as an edge between vertices $i$ and $j$, it is quite natural to regard labeled graphs as monomials which span the invariant ring.
This idea underlies our main result in this section.

\subsection{General algorithm} 

Our first result is the following algorithm, where the input is a certain type of graph with $d$ vertices and $k$ edges, and the output is a basis element in $\P^d_k(\Psi)^G$.
The algorithm takes the same form for all three classical groups.
(See also the appendix, where we include Mathematica code to implement the algorithm.)
Recall that to each classical group $G$ we have already associated (in Theorems~\ref{FTs GLn}--\ref{FTs Sp2n}) a polynomial ring $R$ and the quadratics $r_{ij} \in R$.
Immediately after stating the general algorithm, we will specify the set $\mathcal{G}^d_k$ of graphs for each group.
The key to the algorithm is a $G$-equivariant linear operator $\varphi:R \longrightarrow \P^d(\Psi)$, which (following \cite[p.~174]{sturmfels}) we call the \emph{umbral operator}.
Again, immediately below the algorithm, we explicitly define $\varphi$ for each group.

\begin{tcolorbox}
\begin{alg}[for $G = \GL_n$, $\O_n$, or $\Sp_{2n}$]\label{alg}\ \\

\textbf{Input:} Natural numbers $n,d,k$.

\textbf{Output:} Inside the stable range~\eqref{stable range}, a linear basis for $\P^d_k(\Psi)^{G}$; otherwise, a spanning set.

\bigskip

\begin{enumerate}
\setlength{\itemsep}{10pt}
    \item Let $\mathcal G^d_k$ be the set of graphs (of type determined by $G$) with $d$ vertices and $k$ edges.    
    For each graph $\Gamma \in \mathcal G^d_k$, choose a vertex labeling $1, \ldots, d$ and let $A^\Gamma$ be its adjacency matrix, where the $i$th row and column correspond to the vertex labeled $i$. 
    
    \item For each $\Gamma \in \mathcal G^d_k$, expand $s(\Gamma) \coloneqq\prod_{i,j} r_{ij}^{A^\Gamma_{ij}}$, with $r_{ij}$ as defined in Theorems~\ref{FTs GLn}--\ref{FTs Sp2n}.
    \item The set $\left\{\varphi \circ s(\Gamma) \: \big| \: \Gamma \in \mathcal G^d_k\right\}$ is a spanning set (in the stable range, a basis) for $\P^d_k(\Psi)^G$.
\end{enumerate}
\end{alg}
\end{tcolorbox}

\begin{remark}
In Step 2, the expression $s(\Gamma)$ is the  \emph{symbolic notation} for a $G$-invariant function, in the sense of 19th-century classical invariant theory; hence our ``$s$'' notation.  See \cite[p.~173]{sturmfels} or \cite[Ch.~6]{olver} for discussions of the symbolic method.  In Step 3, the umbral operator $\varphi$ then translates the symbolic notation for an invariant into an explicit expression in terms of coefficients.
\end{remark}

\begin{example}
    To gain some intuition for Algorithm~\ref{alg} applied to an individual graph, let $G = \GL_n$, and set the parameters $n=4$, $d = 3$, and $k = 5$:

    \begin{enumerate}
        \item Start with the following digraph $\Gamma \in \mathcal{G}^3_5$ (defined below in Section~\ref{sec:GLalg}); choosing to label the vertices $1, 2, 3$ from left to right, we then obtain the adjacency matrix $A^\Gamma$:
        \[
        \Gamma = \;\begin{tikzpicture}[>=stealth, scale=0.8, baseline]
\GraphInit[vstyle=Simple]
\tikzset{VertexStyle/.append style={minimum size=6pt, inner sep=1pt}}
\Vertex{A} \EA(A){B} \EA(B){C}
\Loop[dist=15pt, style={thick,->}, dir=EA](C)
\Edge[style={thick,->, bend left}](A)(B)
\Edge[style={thick,->, bend left}](B)(C)
\Edge[style={thick,->, bend left=50}](B)(C)
\Edge[style={thick,->, bend left}](C)(B)
\end{tikzpicture} \qquad \leadsto \qquad A^\Gamma = 
        \begin{bmatrix}
        0 & 1 & 0 \\
        0 & 0 & 2 \\
        0 & 1 & 1
        \end{bmatrix}.\]

    \item Reading $A^\Gamma$ as the degree matrix of a monomial in the quadratics $r_{ij}$, we obtain
    \[
        s(\Gamma) = r_{12} r_{23}^2 r_{32} r_{33}.
    \]
    
    \item Following~\eqref{rij GLn}, the expansion of $s(\Gamma)$ contains 640 terms in the variables $x_{ij}$ and $y_{ij}$.
    For example, one of these terms is $x_{22} x_{32}x_{43}^3 y_{21} y_{42}^2 y_{33} y_{43}$.
    The umbral operator $\varphi$, as defined below in~\eqref{phiGLn}, transforms this term into a product of $\chat$'s in $\P^3_5(\Psi)^{\GL_4}$:
    \[
    \varphi: x_{22} x_{32}x_{43}^3 y_{21} y_{42}^2 y_{33} y_{43} \longmapsto \chat_{(0000),(0100)} \chat_{(0110),(0002)} \chat_{(0003),(0011)}.
    \]
    Recalling the definition of $\chat_{\al,\be}$ from Section~\ref{sec:notation}, we can spell out this term concretely as a function on $\Psi = \C[x_1, \ldots, x_4, \d_1, \ldots, \d_4]$, as follows:
    \[
    6 \cdot (\text{coeff.~of $\d_2$}) \cdot (\text{coeff.~of $x_2x_3\d_4^2$}) \cdot (\text{coeff.~of $x_4^3\d_3 \d_4$}).
    \]
    Carrying this out for all terms in $s(\Gamma)$ and taking the sum, one obtains the invariant $\varphi \circ s(\Gamma)$.
    Since $n \geq d$ in this example, we are in the stable range~\eqref{stable range}, and so $\Gamma \mapsto \varphi \circ s(\Gamma)$ is a bijection between $\mathcal{G}^3_5$ and a linear basis for $\P^3_5(\Psi)^{\GL_4}$.
    \end{enumerate}
\end{example}

\subsection{Details: the general linear group}\label{sec:GLalg}

For $G = \GL_n$, we define $\mathcal G^d_k$ to be the set of \emph{directed} multigraphs with loops, with $d$ vertices and $k$ edges, up to isomorphism.
The quadratics $r_{ij}$ are defined in~\eqref{rij GLn}.
Define the \emph{umbral operator} $\varphi: R \longrightarrow \P^d(\Psi)$ by
\begin{equation}\label{phiGLn}
    \varphi: \prod_{i,j} x_{ij}^{p_{ij}}y_{ij}^{q_{ij}} \longmapsto \prod_{j=1}^d \chat_{p_{_{\bullet j,}}\: q_{_{\bullet j}}}
\end{equation}
where $p_{_{\bullet j}} = \left(p_{1j},\ldots,p_{nj}\right)$ and $q_{_{\bullet j}} = \left(q_{1j},\ldots,q_{nj}\right)$, and extend by linearity.

\subsection{Details: the orthogonal group}

For $G = \O_n$, we define $\mathcal G^d_k$ to be the set of \emph{undirected} multigraphs with loops, with $d$ vertices and $k$ edges, up to isomorphism.
The quadratics $r_{ij}$ are defined in~\eqref{rij On}.
Define the \emph{umbral operator} $\varphi: R \longrightarrow \P^d(\Psi)$ by
\begin{equation}\label{phiOn}
    \varphi:\prod_{i,j} x_{ij}^{p_{ij}} \longmapsto \prod_{j=1}^d \chat_{p_{_{\bullet j}}}
\end{equation}
where $p_{\bullet j} = \left(p_{1j},\ldots,p_{nj}\right)$, and extend by linearity.

\subsection{Details: the symplectic group}\label{section:Sp2nAlg}

For $G = \Sp_{2n}$, the graphs in $\mathcal G^d_k$ are more delicate to define.  
Let $\Gamma$ be an undirected multigraph with no loops, having $d$ vertices.  
Arbitrarily choose vertex labels $1, \ldots, d$; then each element $\sigma$ of the symmetric group $\Sd$ determines another labeled graph $\sigma \cdot \Gamma$, obtained by permuting the original labels. 
The stabilizer $\operatorname{stab}_{\Sd}(\Gamma)$ contains all elements $\sigma \in \Sd$ such that $\Gamma$ is isomorphic to $\sigma \cdot \Gamma$ as a labeled graph.
 If $\{i,j\}$ is an edge of $\Gamma$ with $i < j$, then we say that $\sigma$ \emph{inverts} $\{i,j\}$ 
 if $\sigma(i) > \sigma(j)$.  
 We now define the following set of graphs for the case $G=\Sp_{2n}$:
 \begin{equation}
     \label{GdkSp2n}
     \mathcal G^d_k \coloneqq \left\{ \parbox{5.5cm}{\centering $\Gamma$ a loopless undirected multigraph\\ with $d$ vertices and $k$ edges} \: \middle| \: \parbox{5.5cm}{\centering for all $\sigma \in \operatorname{stab}_{\Sd}(\Gamma)$,\\$\sigma$ inverts an even number of edges} \right\}.
 \end{equation}
 This set is independent of the initial choice of labeling for $\Gamma$, since $\Sd$ acts $d$-transitively on the vertex labels.  
 Below we present examples of the graphs in $\mathcal G^d_k$ for $d \leq 3$.  
 (For $d=1$, by definition, $\mathcal G^1_k = \varnothing$ since loops are not allowed.)

\begin{example}

Let $d=2$, so that any loopless graph $\Gamma$ consists of $k$ edges connecting the two vertices.  
Then both elements of $\mathfrak S_2$ stabilize the labeled graph $\Gamma$. 
Since the nontrivial element of $\mathfrak S_2$ inverts all $k$ edges, we have that $\Gamma \in \mathcal G^2_k$ if and only if $k$ is even.  
Hence $\# \mathcal G^2_k$ equals $1$ if $k$ is even, and $0$ if $k$ is odd.
\end{example}

\begin{example}
Now suppose $d=3$.  Below we list the elements of $\mathcal G^3_k$ for the first few values of $k$.  The reader can check that we have included precisely the graphs defined in \eqref{GdkSp2n}:
\begin{align*}
    \mathcal G^3_0 &= \left\{ \: \v{\begin{tikzpicture}[scale=0.4, rotate=90]\GraphInit[vstyle=Simple]
    \tikzset{VertexStyle/.append style={minimum size=6pt, inner sep=1pt}}
    \Vertices{circle}{A,B,C} \end{tikzpicture}}\: \right\}.\\
\mathcal G^3_1 &= \varnothing.\\
\mathcal G^3_2 &= \left\{
    \v{
    \begin{tikzpicture}[>=stealth, scale=0.4, rotate=-30]
\GraphInit[vstyle=Simple]
\tikzset{VertexStyle/.append style={minimum size=6pt, inner sep=1pt}}
\Vertices{circle}{A,B,C}
\tikzset{EdgeStyle/.style = {-, very thick}}
\Edge(A)(B) \Edge(B)(C)
\end{tikzpicture}
}, \quad  \v{
    \begin{tikzpicture}[>=stealth, rotate=-150, scale=0.4]
\GraphInit[vstyle=Simple]
\tikzset{VertexStyle/.append style={minimum size=6pt, inner sep=1pt}}
\Vertices{circle}{A,B,C}
\tikzset{EdgeStyle/.style = {-, very thick, bend left}}
\Edge(A)(B) \Edge(B)(A)
\end{tikzpicture}
}\right\}.\\
\mathcal G^3_3 &= \left\{
\v{
    \begin{tikzpicture}[>=stealth]
\GraphInit[vstyle=Simple]
\tikzset{VertexStyle/.append style={minimum size=6pt, inner sep=1pt}}
\Vertex{A} \EA(A){B} \EA(B){C}
\tikzset{EdgeStyle/.style = {-, very thick}}
\Edge(B)(C)
\tikzset{EdgeStyle/.style = {-, very thick, bend left}}
\Edge(A)(B) \Edge(B)(A)
\end{tikzpicture}
}\right\}.\\
\mathcal G^3_4 &= \left\{ 
\v{
    \begin{tikzpicture}[>=stealth]
\GraphInit[vstyle=Simple]
\tikzset{VertexStyle/.append style={minimum size=6pt, inner sep=1pt}}
\Vertex{A} \EA(A){B} \EA(B){C}
\tikzset{EdgeStyle/.style = {-, very thick}}
\Edge(B)(C) \Edge(A)(B)
\tikzset{EdgeStyle/.style = {-, very thick, bend left}}
\Edge(A)(B) \Edge(B)(A)
\end{tikzpicture}
}, \quad \v{
    \begin{tikzpicture}[>=stealth]
\GraphInit[vstyle=Simple]
\tikzset{VertexStyle/.append style={minimum size=6pt, inner sep=1pt}}
\Vertex{A} \EA(A){B} \EA(B){C}
\tikzset{EdgeStyle/.style = {-, very thick, bend left}}
\Edge(A)(B) \Edge(B)(A) \Edge(B)(C) \Edge(C)(B)
\end{tikzpicture}
}, \quad \v{
    \begin{tikzpicture}[>=stealth, rotate=-150, scale=0.4]
\GraphInit[vstyle=Simple]
\tikzset{VertexStyle/.append style={minimum size=6pt, inner sep=1pt}}
\Vertices{circle}{A,B,C}
\tikzset{EdgeStyle/.style = {-, very thick, bend left}}
\Edge(A)(B) \Edge(B)(A)
\tikzset{EdgeStyle/.style = {-, very thick, bend left = 10}}
\Edge(A)(B) \Edge(B)(A)
\end{tikzpicture}
}, \quad
\v{
    \begin{tikzpicture}[>=stealth, rotate=-150, scale=0.4]
\GraphInit[vstyle=Simple]
\tikzset{VertexStyle/.append style={minimum size=6pt, inner sep=1pt}}
\Vertices{circle}{A,B,C}
\tikzset{EdgeStyle/.style = {-, very thick}}
\Edge(A)(C) \Edge(C)(B)
\tikzset{EdgeStyle/.style = {-, very thick, bend left}}
\Edge(A)(B) \Edge(B)(A)
\end{tikzpicture}
}\right\}.\\
     \mathcal G^3_5 & = \left\{ \v{
    \begin{tikzpicture}[>=stealth]
\GraphInit[vstyle=Simple]
\tikzset{VertexStyle/.append style={minimum size=6pt, inner sep=1pt}}
\Vertex{A} \EA(A){B} \EA(B){C}
\tikzset{EdgeStyle/.style = {-, very thick}}
\Edge(A)(B)
\tikzset{EdgeStyle/.style = {-, very thick, bend left}}
\Edge(A)(B) \Edge(B)(A) \Edge(B)(C) \Edge(C)(B)
\end{tikzpicture}
},\quad \v{
    \begin{tikzpicture}[>=stealth]
\GraphInit[vstyle=Simple]
\tikzset{VertexStyle/.append style={minimum size=6pt, inner sep=1pt}}
\Vertex{A} \EA(A){B} \EA(B){C}
\tikzset{EdgeStyle/.style = {-, very thick}}
\Edge(B)(C)
\tikzset{EdgeStyle/.style = {-, very thick, bend left}}
\Edge(A)(B) \Edge(B)(A) 
\tikzset{EdgeStyle/.style = {-, very thick, bend left = 10}}
\Edge(A)(B) \Edge(B)(A)
\end{tikzpicture}
}\right\}.
\end{align*}
\end{example}
As for the remaining details in Algorithm \ref{alg} for $G=\Sp_{2n}$: the quadratics $r_{ij}$ are defined in~\eqref{rij Sp2n}, and we define the \emph{umbral operator} $\varphi: R \longrightarrow \P^d(\Psi)$ by
\begin{equation}\label{phiSp2n}
    \varphi:\prod_{i,j} x_{ij}^{p_{ij}} \longmapsto \prod_{j=1}^d \chat_{p_{_{\bullet j}}}
\end{equation}
where $p_{\bullet j} = \left(p_{1j},\ldots,p_{2n,j}\right)$, and extend by linearity.

\section{Proof of Algorithm~\ref{alg}}

The essence of the proof is the symmetrization of Weyl's fundamental theorems: that is, we restrict the isomorphism $\mu^\sharp$ (as defined in Theorems~\ref{FTs GLn}--\ref{FTs Sp2n}) to the subspace of invariants under the action of the symmetric group $\Sd$.
Following~\cite[p.~174]{sturmfels}, on any $\Sd$-module we denote the \emph{symmetrization} operator (i.e., the Reynolds operator for $\Sd$) by a star, so that
\[
x^* \coloneqq \frac{1}{d!}\sum_{\sigma \in \Sd} \sigma \cdot x.
\]

We first observe that $\mu^\sharp$ is $\Sd$-equivariant, as follows.
On one hand, $\Sd$ acts on $\M_d$ by simultaneous permutations of rows and columns (i.e., via conjugation by permutation matrices).
This makes $\P(\M_d)$ into an $\Sd$-module, where $\sigma \cdot z_{ij} = z_{\sigma^{-1}(i), \sigma^{-1}(j)}$ for $\sigma \in \Sd$.
On the other hand, $\Sd$ acts on $V^{\oplus d}$ by permuting copies of $V$, and this makes the polynomial ring $R$ into an $\Sd$-module, via $\sigma \cdot x_{ij} = x_{i,\sigma^{-1}(j)}$ and $\sigma \cdot y_{ij} = y_{i, \sigma^{-1}(j)}$.
It is easy to check that $\mu^\sharp$ intertwines these two $\Sd$-actions.

We provide full details for $\GL_n$, and then simply note the adjustments required for $\O_n$ and $\Sp_{2n}$.
Experts will recognize several of our isomorphisms as examples of \emph{polarization} and \emph{restitution}, to use the language of classical invariant theory~\cite[\S1.2]{dolgachev}.

\begin{proofGL}

Let $R_k \subset R$ denote the component consisting of polynomials which are bihomogeneous of degree $k$ in the $x_{ij}$ and degree $k$ in the $y_{ij}$.
Note that $\mu^\sharp: \P(\M_d^{\leq n}) \longrightarrow R^G$ restricts to a linear isomorphism between graded components $\P^k(\M_d^{\leq n}) \longrightarrow R_k^G$.
Now we consider the following diagram, where we claim that all six vertical arrows (in particular, the two thick arrows on the right-hand side) are linear isomorphisms:
\begin{equation}
    \label{diagram GLn}
\begin{tikzcd}[row sep = 20pt, column sep = 10pt]
    &[10pt] \P(\M_d^{\leq n}) \ar[r, phantom, "\supset"] \ar[d, "\mu^\sharp"'] & \P(\M_d^{\leq n})^{\Sd} \ar[r, phantom, "\supset"] \ar[dd] & \P^k(\M_d^{\leq n})^{\Sd} \ar[dd, very thick] \\
     & R^G \ar[rd, dash, lightgray, thick] \ar[rd, phantom, "\supset" rotate=-35] & & \\
    R \ar[ddr, "\varphi"', bend right] \ar[ru, dash, lightgray, thick] \ar[ru, phantom, "\supset" rotate=35] \ar[rd, dash, lightgray, thick] \ar[rd, phantom, "\supset" rotate=-35] & & R^{G \times \Sd} \ar[r, phantom, "\supset"] \ar[dd] & R^{G \times \Sd}_k \ar[dd, very thick] &\\
    & R^{\Sd} \ar[ru, dash, lightgray, thick] \ar[ru, phantom, "\supset" rotate=35] \ar[d, "\widetilde\varphi"'] & &  \\
    & \P^d(\Psi) \ar[r, phantom, "\supset"] & \P^d(\Psi)^G \ar[r, phantom, "\supset"] & \P^d_k(\Psi)^G &
\end{tikzcd}
\end{equation}
The two vertical arrows to the right of the map $\mu^\sharp$ are successive restrictions of $\mu^\sharp$.
Since $\mu^\sharp$ is $\Sd$-equivariant, its first restriction is a linear isomorphism; moreover, as observed above the diagram, the second restriction of $\mu^\sharp$ is also a linear isomorphism between graded components.

The three vertical arrows to the right of the map $\varphi$ in~\eqref{diagram GLn} are successive restrictions of $\varphi$.
It is clear from its definition in \eqref{phiGLn} that $\varphi:R \longrightarrow \P^d(\Psi)$ is $\Sd$-invariant.  
In fact, its restriction to the subspace of $\Sd$-invariants defines a linear isomorphism $\widetilde \varphi : R^{\Sd} \longrightarrow \P^d(\Psi)$,
which we can see by exhibiting its inverse: for $\al^1,\be^1,\ldots,\al^d,\be^d \in \N^n$, we have
\[
\widetilde \varphi^{-1}: \prod_{j=1}^d \chat_{\al^j,\be^j} \longmapsto \left(\prod_{j=1}^d \prod_{i=1}^n x_{ij}^{\al^j_i} y_{ij}^{\be^j_i}\right)^*.
\]
Next, in order to prove that this isomorphism restricts to an isomorphism between $G$-invariant subspaces, it suffices to show that $\varphi$ is $G$-equivariant.  
As a $G$-module, we have
\[
R \cong \P\!\left(V^{\oplus d} \oplus (V^*)^{\oplus d}\right) \cong \textstyle\bigotimes^d \P(V \oplus V^*) = \bigotimes^d \Psi.
\]
This equivalence of $G$-modules $R \longrightarrow \bigotimes^d \Psi$ is given by
\[
\prod_{j=1}^d \prod_{i=1}^n  x_{ij}^{p_{ij}}y_{ij}^{q_{ij}} \longmapsto \bigotimes_{j=1}^d \prod_{i=1}^n x_i^{p_{ij}}\d_i^{q_{ij}}
\]
on monomials, and extended by linearity.
Identifying $\Psi$ with $\Psi^\circ$ as in \eqref{dualbasis}, we now have an equivalence of $G$-modules defined by the isomorphism
\begin{align*}
    R &\longrightarrow \textstyle\bigotimes^d \Psi^\circ,\\
    \prod_{i,j} x_{ij}^{p_{ij}}y_{ij}^{q_{ij}} &\longmapsto \bigotimes_{j=1}^d \chat_{p_{_\bullet j,}\:q_{_\bullet j}}.
\end{align*}
Composing this with the canonical projection $\bigotimes^d \Psi^\circ \longrightarrow S^d(\Psi^\circ) \cong \P^d(\Psi)$ preserves $G$-equivariance, and in fact yields the umbral operator $\varphi$.
Hence $\varphi$ is $G$-equivariant, and so the middle arrow in the bottom row of~\eqref{diagram GLn} is an isomorphism.
Finally, it is clear from~\eqref{phiGLn} and Definition~\ref{def:weight} that $\varphi$ carries monomials in $R_k$ to monomials of weight $k$, that is, monomials in the component $\P^d_k$.
Hence the second thick arrow in~\eqref{diagram GLn} is a linear isomorphism.

To bring the graphs into the picture, we consider the following diagram, where the three steps of Algorithm~\ref{alg} are denoted by the blue circled numbers:
\begin{equation}
    \label{diagram graph GLn}
\begin{tikzcd}
    \widehat{\mathcal{G}}^d_k \ar[d, "z"'] \ar[dd, white, "s"' near start, blue!50, bend right=75] \ar[dd, blue!50, "2"' draw, circle, outer sep = 5pt, midway, bend right=75] \ar[r, twoheadrightarrow, shift right] & \mathcal{G}^d_k \ar[l,blue!50, "1"' draw, circle, outer sep=5pt, shift right] \ar[d, dashed] \\
    \P^k(\M_d^{\leq n}) \ar[r, "*", twoheadrightarrow] \ar[d, "\mu^\sharp"'] & \P^k(\M_d^{\leq n})^{\Sd} \ar[d, very thick] \\
    R^G_k \ar[r, "*", twoheadrightarrow] \ar[rd, white, "\varphi" blue!50, bend right] \ar[rd, blue!50, "3"' draw, circle, outer sep = 5pt, bend right] & R^{G \times \Sd}_k \ar[d, very thick]\\
    & \P^d_k(\Psi)^G    
\end{tikzcd}
\end{equation}
(By slight abuse of notation, we write $\mu^\sharp$ and $\varphi$ to denote their restrictions above.)
Let $\widehat{\mathcal{G}}^d_k$ denote the set of \emph{labeled} digraphs, with vertices labeled $1, \ldots, d$.
Letting $A^\Gamma$ denote the adjacency matrix of a graph $\Gamma \in \widehat{\mathcal{G}}^d_k$, we define
\[
z(\Gamma) \coloneqq \prod_{i,j=1}^d z_{ij}^{A^\Gamma_{ij}}.
\]
If $n \geq d$, then $\M_d^{\leq n} = \M_d$ and so $z(\widehat{\mathcal{G}}^d_k)$ is a basis for $\P^k(\M_d^{\leq n})$.
In general, $\P(\M^{\leq n}_d)$ has a basis consisting of monomials of width $\leq n$, with ``width'' in the sense of~\cite[Lemma 7]{SturmfelsGB}; but the width of a monomial is necessarily less than or equal to its degree, and so if $n \geq k$, then $z(\widehat{\mathcal{G}}^d_k)$ is still a basis for $\P^k(\M^{\leq n}_d)$.  
Therefore, in the stable range $n \geq \min\{d,k\}$, the set $z(\widehat{\mathcal{G}}^d_k)$ is a basis for $\P^k(\M^{\leq n}_d)$; otherwise, it is only a spanning set.

Note that $\Sd$ acts on $\widehat{\mathcal{G}}^d_k$ by permuting vertex labels, and that $z$ is $\Sd$-equivariant.
Hence the symmetrization $z(\Gamma)^* \in \P^k(\M_d^{\leq n})^{\Sd}$ is independent of the labeling of $\Gamma$, and the entire diagram commutes.
In particular, $z(\Gamma)^*$ is well defined on $\mathcal{G}^d_k$ (this is the dashed arrow in the diagram), and its image is a basis for $\P^k(\M_d^{\leq n})^\Sd$ inside the stable range.
Composing with the restriction of $\mu^\sharp$ followed by the restriction of $\varphi$ (i.e., the isomorphisms depicted by the two thick arrows), we obtain the desired basis for $\P^d_k(\Psi)^G$.
By the commutativity of the diagram, the sequence of three vertical arrows on the right-hand side of~\eqref{diagram graph GLn} is equivalent to the three steps of Algorithm~\ref{alg} (shown by the circled numbers 1,2,3).
\end{proofGL}

\begin{proofO}

The entire proof for $G = \GL_n$ goes through, \emph{mutatis mutandis}.
For $r_{ij}$, $\mu^\sharp$, and $\varphi$, we use the definitions in~\eqref{rij On}, \eqref{muOn}, and $\eqref{phiOn}$, respectively.
In the diagram~\eqref{diagram GLn}, we replace $\M_d^{\leq n}$ by $\SM_d^{\leq n}$, and then in all the displayed equations we simply remove the $y_{ij}$'s, and replace each $\chat_{\alpha,\beta}$ by $\chat_\alpha$.
In the diagram~\eqref{diagram graph GLn} we treat $\mathcal{G}^d_k$ as a set of \emph{un}directed graphs, and we restrict the map $z$ to just the upper-triangular entries, via
\[
z(\Gamma) \coloneqq \prod_{1 \leq i \leq j \leq d} z_{ij}^{A^\Gamma_{ij}} \in \P^k(\SM_d^{\leq n}).
\]
The rest of the proof proceeds identically.
\end{proofO}

\begin{proofSp}

Once again, this is just a matter of changing the obvious details in the proof for $G = \GL_n$.
For $r_{ij}$, $\mu^\sharp$, and $\varphi$, we use the definitions in~\eqref{rij Sp2n}, \eqref{muSp2n}, and \eqref{phiSp2n}, respectively.
In the diagram~\eqref{diagram GLn} we replace $\M_d^{\leq n}$ by $\AM_d^{\leq 2n}$.
Note that the stable range is different here, namely $n \geq \min\{ \lfloor d/2 \rfloor, \lfloor k/2 \rfloor\}$, because the matrices in the determinantal variety have rank $\leq 2n$.
Recall that~\eqref{GdkSp2n} defines the set $\mathcal{G}^d_k$.
For a labeled graph $\Gamma \in \widehat{\mathcal{G}}^d_k$, let $A^{\Gamma} \in \AM_d(\Z)$ be the skew-symmetric matrix whose upper-triangular entries are those of the adjacency matrix of $\Gamma$.  
We define the map $z$ by
\[
z(\Gamma) \coloneqq \prod_{1 \leq i < j \leq d} z_{ij}^{A^{\Gamma}_{ij}} \in \P^k(\AM^{\leq 2n}_d).
\]

The delicate point in this proof is to show that our definition~\eqref{GdkSp2n} determines the correct subset $\mathcal{G}^d_k$ of graphs in order to obtain bases.
Let $\Gamma$ be a labeled loopless graph with $d$ vertices and $k$ edges.
Due to the relation $z_{ij}=-z_{ji}$ in $\P^k(\AM^{\leq 2n}_d)$, it follows that
\begin{equation}\label{extrasign}
\sigma \cdot z(\Gamma) = (-1)^{\operatorname{inv}(\Gamma,\sigma)}\: z(\sigma \cdot \Gamma),
\end{equation}
where $\operatorname{inv}(\Gamma,\sigma)$ denotes the number of edges of $\Gamma$ inverted by $\sigma$.  
Unlike the cases for $\GL_n$ and $\O_n$, these signs make it possible for the symmetrization operator to kill $z(\Gamma)$ for certain graphs $\Gamma$.

It is true that a spanning set for $\P^k(\AM_d^{\leq 2n})$ is certainly given by the set of all $z(\Gamma)$ where $\Gamma$ ranges over \textit{all} loopless graphs with $d$ (labeled) vertices and $k$ edges; we now show that our definition of $\mathcal G^d_k$ in \eqref{GdkSp2n} excludes those $\Gamma$ that are killed by the passage to $\Sd$-invariants, i.e., for which $z(\Gamma)^*=0$.  To this end, let $\sigma \in \operatorname{stab}_{\Sd}(\Gamma)$. 
 It follows that $A^{\sigma \cdot \Gamma}=A^{\Gamma}$, and so \begin{equation}\label{ztaugamma}
z(\sigma \cdot \Gamma)=z(\Gamma).
\end{equation}
Thus we have
\begin{align*}
    z(\Gamma)^* &\coloneqq \frac{1}{d!} \sum_{\tau \in \Sd} \tau \cdot z(\Gamma) \\
    &= \frac{1}{d!}\sum_{\tau \in \Sd} \tau \cdot z(\sigma \cdot \Gamma) && \text{by \eqref{ztaugamma}}\\
    &= \frac{1}{d!}\sum_{\tau \in \Sd} \tau \cdot [(-1)^{{\rm inv}(\Gamma,\sigma)} \sigma \cdot z(\Gamma)] && \text{by \eqref{extrasign}}\\
    &= (-1)^{{\rm inv}(\Gamma,\sigma)} \cdot \frac{1}{d!} \sum_{\tau \in \Sd} \tau\sigma \cdot z(\Gamma) \\
    &= (-1)^{{\rm inv}(\Gamma,\sigma)} \cdot \frac{1}{d!} \sum_{\tau \in \Sd} \tau \cdot z(\Gamma) \\
    &= (-1)^{{\rm inv}(\Gamma,\sigma)} \cdot z(\Gamma)^*.
\end{align*}
Therefore $z(\Gamma)^* = 0$ if and only if ${\rm inv}(\Gamma, \sigma)$ is odd for some $\sigma \in {\rm stab}_{\Sd}(\Gamma)$.
This justifies our definition~\eqref{GdkSp2n} of the set $\mathcal{G}^d_k$ for the symplectic group.    
\end{proofSp}

\section{Hilbert series via branching multiplicities}
\label{sec:branching}

\subsection{Combinatorial difficulty outside the stable range}

Inside the stable range~\eqref{stable range}, one could scarcely hope for a nicer combinatorial understanding of $\P(\Psi)^G$ than that given by Algorithm~\ref{alg}:
specifically, we can obtain the dimension of each graded component $\P^d_k(\Psi)^G$ simply by counting the graphs in $\mathcal{G}^d_k$.
\textit{Outside} the stable range, however, the graph-counting approach overshoots the true dimension, due to the complicated linear dependencies arising among the images $\varphi \circ s(\Gamma)$.

Taking $G = \GL_n$, for example, it is straightforward to give a linear basis for $\P^k(\M_d^{\leq n})$ consisting of the degree-$k$ monomials of \emph{width} at most $n$; the term ``width'' follows the sense of Sturmfels~\cite{SturmfelsGB}, which we now summarize.
A pair $(T,U)$ of semistandard Young tableaux with the same shape, each having (say) $\ell$ columns, and filled with entries from the set $[d] \coloneqq \{1, \ldots, d\}$, can be viewed as a product of $\ell$ determinants in $\P(\M_d)$, as follows: the $i$th determinant is the minor of matrix coordinates whose rows (resp., columns) are given by the $i$th column of $T$ (resp., $U$).
It follows that the size of $T$ (equivalently, of $U$) is the degree of the corresponding function in $\P(\M_d)$.
The set of all such pairs (often called \emph{bitableaux} in the literature) thus yields a linear basis for $\P(\M_d)$.
Via the straightening law of~\cite{DRS}, it can be shown that by restricting to those bitableaux with size $k$ and at most $n$ rows, one obtains a linear basis for $\P^k(\M^{\leq n}_d)$.
Sturmfels showed that these somewhat complicated basis elements (an example of \emph{standard monomials}) can be replaced by certain ordinary monomials in the matrix coordinates $z_{ij}$, by applying the Robinson--Schensted--Knuth (RSK) correspondence.
The RSK correspondence, described by Knuth in~\cite[\S3]{Knuth}, is a bijection between the set of bitableaux and the set $\M_d(\N)$.
In particular, the sum of the entries in the matrix ${\rm RSK}(T,U)$ equals the size of $T$; moreover, when the support of  ${\rm RSK}(T,U)$ is viewed as a subset of the poset $[d] \times [d]$, the \emph{width} of that support (i.e., the size of the largest antichain) equals the number of rows in $T$.
By viewing ${\rm RSK}(T,U)$ as the degree matrix of a monomial in the variables $z_{ij}$, we can thus associate each bitableau to an ordinary monomial rather than to a standard monomial.
We show this below in an example where $d=5$ and $n \geq 3$:
\[
\left(
\ytableausetup{centertableaux,smalltableaux}
\ytableaushort{11223,335,55}_{\textstyle{,}} \: \ytableaushort{11123,233,44}
\right) \xmapsto{\text{RSK}}   
    \left[\begin{smallmatrix}
        1&\textcolor{gray}{0}&\textcolor{gray}{0}&\textcolor{gray}{0}&2\\
        1&\textcolor{gray}{0}&\textcolor{gray}{0}&\textcolor{gray}{0}&1\\
        \textcolor{gray}{0}&\textcolor{gray}{0}&3&\textcolor{gray}{0}&\textcolor{gray}{0}\\
        \textcolor{gray}{0}&2&\textcolor{gray}{0}&\textcolor{gray}{0}&\textcolor{gray}{0}\\
        \textcolor{gray}{0}&\textcolor{gray}{0}&\textcolor{gray}{0}&\textcolor{gray}{0}&\textcolor{gray}{0}
    \end{smallmatrix}
    \right] 
    \quad \leadsto \quad
    z_{11}z_{15}^2 z_{21}z_{25}z_{33}^3z_{42}.
    \]
In this way, a basis for $\P^k(\M_d^{\leq n})$ is given by those ordinary monomials whose degree matrices have entries summing to $k$ and width at most $n$.
Finally, by viewing the degree matrix as an adjacency matrix, one can represent the basis monomials by digraphs with $k$ edges on $d$ (labeled) vertices.
Variations of the RSK correspondence (see~\cite{Burge} or \cite{Conca94}) can be used to obtain similar graphical bases for $\P^k(\SM_d^{\leq n})$ and $\P^k(\AM_d^{\leq 2n})$, corresponding to the groups $\O_n$ and $\Sp_{2n}$.

The difficulty in the present paper, however, arises from our symmetrization, whereby one ``forgets'' the vertex labels in a graph.
In particular, when $n < d$ and $n < k$, it is no longer clear how to parametrize a basis for $\P^k(\M_d^{\leq n})^{\Sd}$.
After all, the sense of ``width'' defining the labeled graphs is inherited from the poset of matrix coordinates, and symmetrization (by its very nature) nullifies the partial order.
Outside the stable range, therefore (and with the exception of Section \ref{sub:n=1}), we leave it as an open problem to exhibit a combinatorial realization of a basis of $\P^d_k(\Psi)^G$.
(We have had some success in the two border cases $n = d-1$ and $n=d-2$, where it is possible to obtain a basis by deleting those graphs containing certain subgraphs; but the choice of these subgraphs is hardly canonical, and moreover this method seems to be hopeless when $d - n > 2$.)

In this section, we remedy this combinatorial defect by taking a representation-theoretic approach.
Note that our goal can be restated more elegantly as follows:
we wish to understand the bigraded Hilbert series of the invariant ring $\P(\Psi)^G$, namely
\[
H_n(q,t) \coloneqq \sum_{d,k} \left(\dim \P^d_k(\Psi)^G\right) q^d t^k.
\]

\subsection{Branching multiplicities}

The facts in this section are standard in any general reference on representation theory, such as~\cite[\S3.2]{gw}.
A \emph{polynomial representation} of $\GL_d$ is a group homomorphism $\rho: \GL_d \longrightarrow \GL_m$ for some $m$, such that the matrix coordinates $\rho(g)_{ij}$ are polynomials in the entries of $g \in \GL_d$.  (As is typical, we will also use the term ``representation'' for the complex $m$-dimensional vector space on which $\rho(\GL_n)\subset \GL_m$ acts.)  A \emph{partition} is a finite, weakly decreasing sequence of positive integers (\emph{parts}), typically denoted by a lowercase Greek letter.  If $d$ is the sum of the parts of $\la$, then we write $\la \vdash d$.  The number of parts of $\la$ is called its \emph{length}, denoted by $\ell(\la)$.  We write $(d)$ to denote the length-1 partition of size $d$.  We also define the set
$$
\Par(k,m) \coloneqq \Big\{ \la \: \Big | \: \la \vdash k \text{ and } \ell(\la)\leq m\Big\}.
$$
We will speak of the \emph{rows} and  \emph{columns} of a partition, by which we mean the rows and columns of its associated Young diagram.

By the theorem of the highest weight, the irreducible polynomial representations of $\GL_d$ are indexed by the partitions $\la$ such that $\ell(\la)\leq d$.  We write $F^\la_d$ for the irreducible representation of $\GL_d$ with highest weight $\la$.  Now consider the restriction of $\GL_d$ to its subgroup $\Sd$, realized as permutation matrices.  If $\mu \vdash d$, we write $Y^\mu_d$ for the irreducible representation of $\Sd$ which corresponds to $F^\mu_d$ via Schur--Weyl duality.  Upon restriction to $\Sd$, the representation $F^\la_d$ decomposes into a direct sum of irreducible representations $Y^\mu_d$; we write
\begin{equation}\label{blamu}
b^\la_\mu \coloneqq \dim \Hom_{\Sd}\!\left(Y^\mu_d,\:F^\la_d\right)
\end{equation}
to denote the \emph{branching multiplicity} of $Y^\mu_d$ in $F^\la_d$.  Hence we have
\begin{equation}\label{GLdSddecomp}
{\rm Res}^{\GL_d}_{\Sd}F^\la_d \cong \bigoplus_{\mu \vdash d} b^\la_\mu\: Y^\mu_d.
\end{equation}

\begin{theorem}\label{thm:branching}
Let $b^\la_\mu$ be as in \eqref{blamu}.  
Then we have the following, for all values of $n,d,k$:

\begin{enumerate}
    \item Let $G=\GL_n$.  Then $\displaystyle\dim \P^d_k(\Psi)^G = \sum_{\la,\mu} \left(b^\la_\mu\right)_,^2$ where $\la \in \Par(k,\:\min\{d,n\})$ and $\mu \vdash d$.
    
    \item Let $G=\O_n$.  Then $\displaystyle\dim \P^d_k(\Psi)^G = \sum_{\la} b^{\la}_{(d),}$ where $\la \in \Par(2k,\:\min\{d,n\})$ with all row lengths even.
    
    \item Let $G = \Sp_{2n}$.  Then $\displaystyle\dim \P^d_k(\Psi)^G = \sum_\la b^\la_{(d),}$ where $\la \in \Par(2k,\:\min\{d,2n\})$ with all column lengths even.
\end{enumerate}
\end{theorem}

\begin{proofGL}
Recall from the two thick arrows in diagram~\eqref{diagram GLn} that we have a linear isomorphism $\P^k(\M_d^{\leq n})^\Sd \cong \P^d_k(\Psi)^G$.  Hence it will suffice to show that $\dim \P^k(\M_d^{\leq n})^\Sd$ equals the sum in the theorem.

The determinantal variety $\M_d^{\leq n}$ admits an action by $\GL_d \times \GL_d$, via $(g,h)\cdot X = gXh^{-1}$, which extends to $\P(\M_d^{\leq n})$ in the usual way.  
Weyl's second fundamental theorem (SFT, from Section~\ref{sub:FFT SFT}) can be applied to obtain the decomposition below~\cite[Thm.~12.2.12.3]{gw}:
\[
\P^k(\M_d^{\leq n}) \cong \bigoplus_{\la\in \Par(k,\:\min\{d,n\})} \left(F^\la_d\right)^* \otimes F^\la_d.
\]
Restricting to $\Sd \times \Sd$, we use \eqref{GLdSddecomp} to write
\begin{align*}
\P^k(\M_d^{\leq n}) &\cong \bigoplus_{\la} \left(\bigoplus_{\mu \vdash d} b^\la_\mu \:Y^\mu_d\right)^* \otimes \left(\bigoplus_{\nu \vdash d} b^\la_\nu\: Y^\nu_d\right)\\
& \cong \bigoplus_{\la,\mu,\nu} b^\la_\mu b^\la_\nu \left(Y^\mu_d \otimes Y^\nu_d\right),
\end{align*}
since the irreducible representations of $\Sd$ are self-dual.  Further restricting to the diagonal subgroup $\Delta(\Sd) = \{(\sigma,\sigma) \mid \sigma \in \Sd\}$, we note that the $\Delta(\Sd)$-action is just by conjugation on $\M_d^{\leq n}$.  Therefore we are interested in the subspace of $\Delta(\Sd)$-invariants:
$$
\P^k(\M_d^{\leq n})^\Sd \cong \left[\bigoplus_{\la,\mu,\nu} b^\la_\mu b^\la_\nu  \left(Y^\mu_d \otimes Y^\nu_d\right)\right]^{\Delta(\Sd)} = \bigoplus_{\la,\mu,\nu} b^\la_\mu b^\la_\nu \left( Y^\mu_d \otimes Y^\nu_d\right)^{\Delta(\Sd)}.
$$
But by Schur's lemma, we have $\dim (Y^\mu_d \otimes Y^\nu_d)^{\Delta(\Sd)} = \dim \Hom_{\Sd}(Y^\mu_d,Y^\nu_d) = \delta_{\mu,\nu}$.  Therefore we keep only the summands in which $\mu=\nu$, and we obtain
$$
\P^k(\M_d^{\leq n})^\Sd \cong \bigoplus_{\la,\mu} \left(b^\la_\mu\right)^2 \underbrace{\left(Y^\mu_d \otimes Y^\mu_d\right)^{\Delta(\Sd)}}_{\text{1-dimensional}}
$$
which proves Part 1 of the theorem.
\end{proofGL}

\begin{proofO}
In the $\O_n$-analogue of the diagram \eqref{diagram GLn}, we have an isomorphism $\P^k(\SM_d^{\leq n})^\Sd \cong \P^d_k(\Psi)^G$.  Hence it will suffice to show that $\dim \P^k(\SM_d^{\leq n})^\Sd$ equals the sum in the theorem.

The determinantal variety $\SM_d^{\leq n}$ admits an action by $\GL_n$, via $g \cdot X =  g^{-T}Xg^{-1}$, which extends to $\P(\SM_d^{\leq n})$ in the usual way.  
The SFT can be applied to obtain the decomposition below~\cite[Thm.~ 12.2.14.3]{gw}, which we decompose under the restriction to $\Sd$ using \eqref{GLdSddecomp}:
\begin{align*}
\P^k(\SM_d^{\leq n}) &\cong \bigoplus_{\la} F^{\la}_d\\
&\cong \bigoplus_\la \bigoplus_{\mu \vdash d} b^{\la}_\mu \: Y^\mu_d,
\end{align*}
where the sum ranges over all $\la \in \Par(2k,\: \min\{d,n\})$ with even row lengths.  The partition $(d)$ labels the one-dimensional trivial representation of $\Sd$, and therefore the $\Sd$-invariant subspace above is the sum of the trivial representations $Y^{(d)}_d$ in the sum:
$$
\P^k(\SM_d^{\leq n})^\Sd \cong \bigoplus_\la b^{\la}_{(d)}\:\underbrace{Y^{(d)}_{d.}}_{\mathclap{\text{1-dimensional}}}
$$
This proves Part 2 of the theorem.
\end{proofO}

\begin{proofSp}

As in the previous parts, from the $\Sp_{2n}$-analogue of the diagram \eqref{diagram GLn}, we have an isomorphism $\P^k(\AM_d^{\leq 2n})^\Sd \cong \P^d_k(\Psi)^G$.  Hence we will find the dimension of $\P^k(\AM_d^{\leq 2n})^\Sd$.

The determinantal variety $\AM_d^{\leq 2n}$ admits an action by $\GL_d$, just as in the $\O_n$ case.  
The SFT can be applied to obtain the decomposition below~\cite[Theorem 12.2.15.3]{gw}, which we decompose under the restriction to $\Sd$, using \eqref{GLdSddecomp}:
\begin{align*}
    \P^k(\AM_d^{\leq 2n}) &\cong \bigoplus_{\la} F^\la_d \\
    & \cong \bigoplus_\la \bigoplus_{\mu \vdash d} b^{\la}_\mu \: Y^\mu_{d,}
\end{align*}
where the sum ranges over all $\la \in \Par(2k,\:\min\{d,2n\})$ with even column lengths.  Therefore, as in Part 2, the $\Sd$-invariant subspace decomposes as
$$
\P^k(\AM_d^{\leq 2n})^\Sd \cong \bigoplus_\la b^{\la}_{(d)}\:\underbrace{Y^{(d)}_{d,}}_{\mathclap{\text{1-dimensional}}}
$$
which proves Part 3 of the theorem.
\end{proofSp}

 \begin{remark}
The branching multiplicities $b^\la_\mu$ can be easily programmed by following the algorithm outlined in \cite[\S2]{hsw}.
In \cite[Ex.~7.74]{stanley2}, $b^\la_\mu$ is interpreted as the multiplicity of $F^\la_n$ in $F^\mu(S(V))$, where $F^\mu$ is the Schur functor.  See also \cite{hsw} for an elegant proof using seesaw reciprocity.  It is possible to prove all three parts of Theorem \ref{thm:branching} with this approach, by directly decomposing $\P^d(\Psi)$ and then restricting to the $G$-invariant subspace.
\end{remark}

\subsection{Bases when $\dim V=1$}\label{sub:n=1}  The theme of this section --- and the motivation for looking at branching multiplicities --- has been the extreme difficulty in parametrizing a linear basis for $\P^d_k(\Psi)^G$ outside the stable range.
One exception to this is actually the farthest extreme from the stable range, namely the case when $\dim V = 1$.  
This occurs when $G= \GL_1$ or $\O_1$.  
In this case, the quadratics $r_{ij}$ are monomials, and thus so too are the expressions $s(\Gamma)$ and $\varphi \circ s(\Gamma)$.
It is fairly easy to see from the definitions that when $G = \GL_1$, we have
\[
\varphi \circ s(\Gamma) = \prod_{i = 1}^d \chat_{\operatorname{indeg}(i), \operatorname{outdeg}(i)}
\]
where indeg and outdeg denote the in-degree and out-degree of a vertex of $\Gamma$.  
Likewise, when $G = \O_1$, we have
\[
\varphi \circ s(\Gamma) = \prod_{i=1}^d \chat_{\operatorname{deg}(i)}
\]
where deg denotes the degree of a vertex.  
In other words, for both $\GL_1$ and $\O_1$, the invariant $\varphi \circ s(\Gamma)$ is determined by the degree sequence of $\Gamma$; 
hence we can actually describe a basis for each component $\P^d_k(\Psi)^G$ by replacing the graphs in our method by simpler combinatorial objects.  
For $\GL_1$, we can parametrize a basis by the sets $\{(a_1,b_1), \ldots, (a_d, b_d)\}$ with $a_i,b_i \in \N$ such that $\sum_i a_i = \sum_i b_i = k$.  
For $\O_1$, we can parametrize a basis by the partitions of $2k$ with at most $d$ parts.

For the symplectic group, the parameter $n=1$ corresponds to $\Sp_{2} \cong \SL(2, \C)$.
Since the defining representation $V$ has dimension 2, a combinatorial parametrization of linear bases seems to be much more difficult than for $\GL_1$ and $\O_1$ above.
Nonentheless, thanks to Theorem~\ref{thm:branching}, we do have an especially nice expression for the dimensions in terms of a single branching multiplicity:
\[
\dim \P^d_k(\C[x,y])^{\SL(2, \C)} = b^{(k,k)}_{(d)}
\]
for all $d$ and $k$.
This is because when $n=1$, there is just one element $(k,k)$ in the set of all partitions of $2k$ with length at most 2, and with even column lengths. 

\section{Examples}\label{section:examples}

We have written Mathematica code to implement Algorithm~\ref{alg}, which we used to compute the examples below.
The code is included in the appendix.

\subsection{The general linear group}

For $\GL_n$, we have $\Psi = \P(V\oplus V^*)$, which is isomorphic (as a $\GL_n$-module) to the associated graded algebra $\C[ x_1,\ldots,x_n,\d_1,\ldots,\d_n]$ of the \emph{Weyl algebra}, i.e., the algebra of polynomial-coefficient differential operators on $\P(V)$.  In this context, $x_i$ is the operator that multiplies by $x_i$, and $\d_i$ is the differential operator $\d/\d_i$.  (Of course, the Weyl algebra itself is not commutative, since $\d_i x_i = x_i\d_i+1$, but here we are interested only in the $G$-module structure.)  Explicitly, we have $\operatorname{span}\{x_1,\ldots,x_n\} \cong V^*$ and $\operatorname{span}\{\d_1,\ldots,\d_n\} \cong V$ as $\GL_n$-modules.   (See the detailed exposition in \cite[p.~39]{procesi}.) 
Thus $\Psi$ is of interest because the $\GL_n$-orbits are precisely the equivalence classes of linear partial differential equations with polynomial coefficients, under linear change of coordinates.
In turn, the invariants $\P(\Psi)^{\GL_n}$ are those polynomial functions constant on the $\GL_n$-orbits.
For example ($n=3$), if $\omega(x,y)$ is a polynomial, then the two-dimensional time-dependent Schr\"odinger equation
\[
-i \frac{\d u}{\d t} = \frac{\d^2 u}{\d x^2} + \frac{\d^2 u}{\d y^2} + \omega(x,y)u
\]
corresponds to the Weyl algebra element $i \d_3 + \d_1^2 + \d_2^2 + \omega(x_1, x_2) \in \Psi$, upon setting $x = x_1$, $y = x_2$, and $t = x_3$.
See also Example~\ref{example:vectorfields} below, in which we view vector fields as elements of $\Psi$.

\begin{remark}
One might further hope that the invariants separate the orbits, but this is not true: if the closures of two orbits intersect, then any invariant function is constant on the union of those orbits.
Therefore, the correct statement is the following~\cite[Thm.~3.20]{WallachGIT}: for any classical group $G$, the $G$-invariants separate the \textit{closed} $G$-orbits in $\Psi$, with respect to the following topology on $\Psi$.  We have a filtration $\Psi = \bigcup_{\ell \in \N} \Psi_\ell$, where $\Psi_\ell$ is the space of polynomials of degree at most $\ell$.  A proper subset of $\Psi$ is defined to be closed if it is a Zariski-closed subset of $\Psi_\ell$ for some $\ell$.  The $G$-invariant polynomial functions on $\Psi$ separate the closed $G$-orbits since $G$ does not change degree.  The topological aspects of invariant theory are delicate; see the excellent source \cite{PV}.
\end{remark}

\begin{example}[$d=k=1$]
Let $\Gamma$ be the digraph consisting of a single loop on a single vertex; then we can use Algorithm \ref{alg} to compute the invariant $\varphi \circ s(\Gamma)$ by hand.  Clearly $A^\Gamma = [1]$, and so $s(\Gamma) = r_{11} = \sum_{\ell = 1}^n y_{\ell 1} x_{\ell 1}$.  We have $\varphi(x_{\ell 1} y_{\ell 1}) = \chat_{\epsilon_\ell,\epsilon_\ell}$, where $\epsilon_\ell$ denotes the $n$-tuple whose $\ell$th component is 1 with 0's elsewhere. Summing these up, we obtain the invariant $\sum_{\ell = 1}^n \c_{\epsilon_\ell,\epsilon_\ell}$.  (Note that in this case, the factorials are all $1!$ so that $\chat$ is no different than $\c$.)  To make this more transparent, we write $[\mathbf x^\al \bm{\d}^\be]$ in place of $\c_{\al,\be}$, and obtain the invariant $[x_1\d_1] + \cdots + [x_n\d_n].$
\end{example}

\begin{example}[$n=d=k=2$]\label{example:GL2}
In this example we will compute a basis for $\P^2_2(\Psi)^{\GL_2}$.  Since $G=\GL_2$, we write $x$ and $y$ in place of $x_1$ and $x_2$, so that $\Psi = \C[x,y,\d_x,\d_y]$.  As in the previous example, we write coefficients as $[x^a y^b \d_x ^c \d_y^d]$ rather than $\c_{(a,b),(c,d)}$. 

Note that our parameters $n=d=k=2$ lie in the stable range, and so there is a one-to-one correspondence between $\mathcal G^2_2$ and our desired basis.  There are six directed graphs with two vertices and two edges, and thus $\#\mathcal G^2_2 = 6 = \dim \P^2_2(\Psi)^{\GL_2}$.  In Table \ref{table:GL2} we display each of these graphs, its adjacency matrix (up to simultaneous permutation of rows and columns), and its corresponding basis element in $\P^2_2(\Psi)^{\GL_2}$.
\end{example}

\begin{table}[t]
\begin{tblr}{colspec={|Q[c,m]|Q[c,m]|Q[c,m]|}, rows={25pt}} \hline
$ \Gamma \in \mathcal G^2_2$ & $A^\Gamma$ & $\varphi \circ s(\Gamma)$\\\hline[2pt]
$\v{
    \begin{tikzpicture}[>=stealth, scale=0.8]
\GraphInit[vstyle=Simple]
\tikzset{VertexStyle/.append style={minimum size=6pt, inner sep=1pt}}
\Vertex{A} \EA(A){B}
\Loop[dist=15pt, style={thick,->}](A)
\Loop[dist=15pt,dir=EA, style={thick,->}](B)
\end{tikzpicture}
}$ & $\left[\begin{smallmatrix}1 & 0 \\ 0 & 1 \end{smallmatrix}\right]$ & $[x\d_x]^2 + 2[x\d_x][y\d_y]+[y\d_y]^2$\\ \hline
$\v{
    \begin{tikzpicture}[>=stealth, scale=0.8]
\GraphInit[vstyle=Simple]
\tikzset{VertexStyle/.append style={minimum size=6pt, inner sep=1pt}}
\Vertex{A} \EA(A){B}
\Loop[dist=15pt, style={thick,->}](A)
\Loop[dist=15pt,dir=EA, style={thick,->}](A)
\end{tikzpicture}
}$\hspace{8pt} & $\left[\begin{smallmatrix}2 & 0 \\ 0 & 0 \end{smallmatrix}\right]$ & $[1]\cdot\big(4[x^2\d_x^2]+2[xy\d_x \d_y] + 4[y^2\d_y^2]\big)$\\ \hline
 $\v{
    \begin{tikzpicture}[>=stealth, scale=0.8]
\GraphInit[vstyle=Simple]
\tikzset{VertexStyle/.append style={minimum size=6pt, inner sep=1pt}}
\Vertex{A} \EA(A){B}
\Edge[style={thick,->, bend left}](A)(B)
\Edge[style={thick,->, bend right}](A)(B)
\end{tikzpicture}
}$ & $\left[\begin{smallmatrix}0 & 2 \\ 0 & 0 \end{smallmatrix}\right]$ & $4[x^2][\d_x^2] + 2[xy][\d_x\d_y]+4[y^2][\d_y^2]$\\ \hline
$\v{
    \begin{tikzpicture}[>=stealth, scale=0.8]
\GraphInit[vstyle=Simple]
\tikzset{VertexStyle/.append style={minimum size=6pt, inner sep=1pt}}
\Vertex{A} \EA(A){B}
\Edge[style={thick,->, bend left}](A)(B)
\Edge[style={thick,->, bend left}](B)(A)
\end{tikzpicture}
}$ & $\left[\begin{smallmatrix}0 & 1 \\ 1 & 0 \end{smallmatrix}\right]$ & $[x\d_x]^2 + 2[x\d_y][y\d_x]+[y\d_y]^2$\\ \hline
$\v{
    \begin{tikzpicture}[>=stealth, scale=0.8]
\GraphInit[vstyle=Simple]
\tikzset{VertexStyle/.append style={minimum size=6pt, inner sep=1pt}}
\Vertex{A} \EA(A){B}
\Edge[style={thick,->}](A)(B)
\Loop[dist=15pt, style={thick,->}](A)
\end{tikzpicture}
}$\hspace{8pt} & $\left[\begin{smallmatrix}1 & 1 \\ 0 & 0  \end{smallmatrix}\right]$ & $2[x\d_x^2][x]+[x\d_x\d_y][y]+[y\d_x\d_y][x]+2[y\d_y^2][y]$ \\ \hline
$\v{
    \begin{tikzpicture}[>=stealth, scale=0.8]
\GraphInit[vstyle=Simple]
\tikzset{VertexStyle/.append style={minimum size=6pt, inner sep=1pt}}
\Vertex{A} \EA(A){B}
\Edge[style={thick,->}](B)(A)
\Loop[dist=15pt, style={thick,->}](A)
\end{tikzpicture}
}$\hspace{8pt} & $\left[\begin{smallmatrix}1 & 0 \\ 1 & 0 \end{smallmatrix}\right]$ & $2[x^2\d_x][\d_x]+[xy\d_x][\d_y]+[xy\d_y][\d_x]+2[y^2\d_y][\d_y]$ \\\hline
   \end{tblr}
   \caption{A basis for $\P^2_2(\Psi)^{\GL_2}$, computed using Algorithm \ref{alg}.}
   \label{table:GL2}
   \end{table}

\begin{remark}
We observe a general phenomenon from Table \ref{table:GL2}, which follows directly from the definition of the umbral operator in \eqref{phiGLn}.  Let $\Gamma \in \mathcal G^d_k$, and choose an adjacency matrix $A^\Gamma$ for any vertex labeling.  Then in each term
$$
\c_{\al^1,\be^1} \cdots \c_{\al^d,\be^d} \quad (\al^i,\be^i \in \N^n)
$$
of the corresponding invariant $\varphi \circ s(\Gamma)$, the $d$ factors $\c_{\al^i,\be^i}$ can be rearranged so that
\begin{align}\label{rowcolsums}
    \begin{split}
        |\al^i| &= \text{$i$th column sum of $A^\Gamma$},\\
        |\be^i| &= \text{$i$th row sum of $A^\Gamma$}
    \end{split}
\end{align}
simultaneously for each $i=1,\ldots, d$.  

Now recall that an element $\psi \in \Psi$ is a \emph{vector field} if it takes the form
\[
\psi = \sum_{i=1}^n f_i(\mathbf x) \:\d_i
\]
with each $f_i(\mathbf x) \in \C[\mathbf x]$.  In this case, $\c_{\al,\be}(\psi) \neq 0$ implies that $|\be|=1$.  By the observation \eqref{rowcolsums}, then, $\varphi \circ s(\Gamma)$ vanishes on $\psi$ unless $A^\Gamma$ contains exactly one 1 in each row, with 0's elsewhere; equivalently, every vertex in $\Gamma$ has out-degree 1.  
For fixed $d$, each such digraph can be regarded naturally as the equivalence class of a \emph{finite dynamical system} on a $d$-element set.
\end{remark}

\begin{figure}
    \centering
\begin{subfigure}[b]{0.48\textwidth}
\centering
\includegraphics[width=0.7\linewidth]{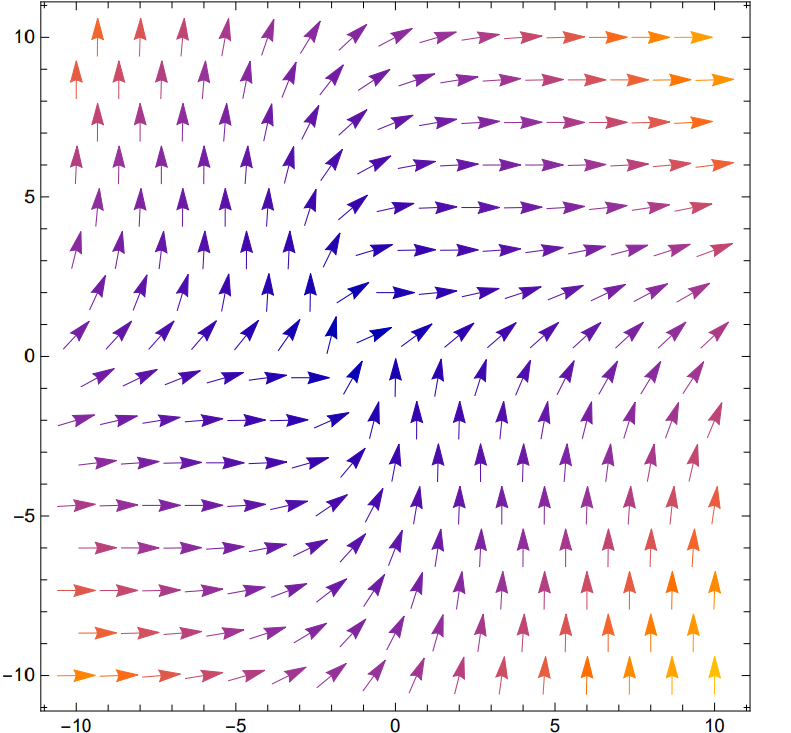} 
\caption{The vector field $\psi$}
\label{subfig:psi}
\end{subfigure}
\begin{subfigure}[b]{0.48\textwidth}
\centering
\includegraphics[width=0.7\linewidth]{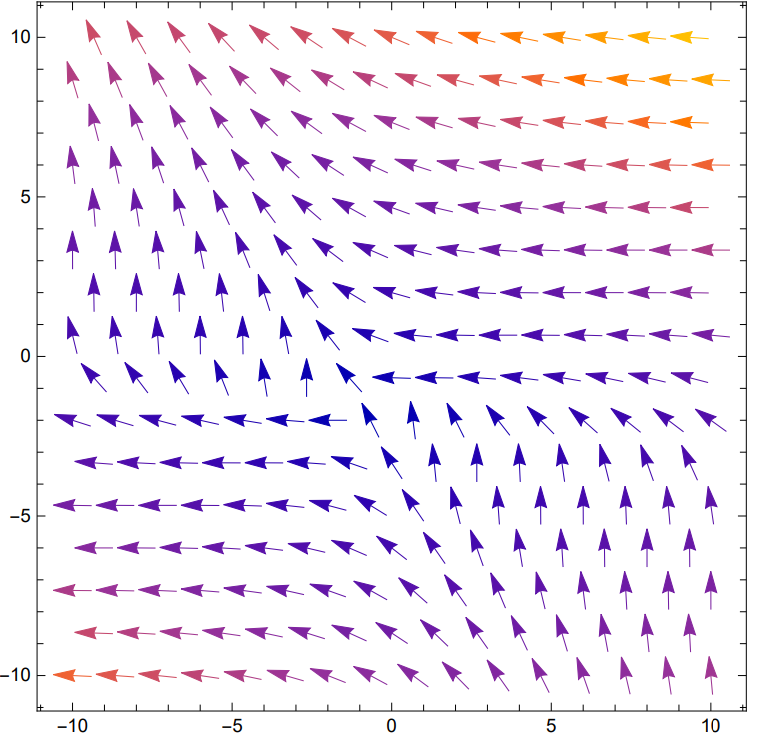} 
\caption{The transformed vector field $g \cdot \psi$}
\label{subfig:gpsi}
\end{subfigure}
    \caption{The vector fields from Example \ref{example:vectorfields}.}
    \label{fig:vectorfields}
\end{figure}

\begin{example}[Quadratic invariants on two-dimensional vector fields]\label{example:vectorfields}

Let $G = \GL_2$, and define the vector field
$$
\psi = \big(x^2 + y^2 + 2xy + 2x + 2y + 1\big)\: \d_x + \big(x^2 + y^2 - 2xy + 4x - 4y+4\big)\: \d_y,
$$
which (since all coefficients are real) we depict in Figure \ref{subfig:psi}.  Let $g = \left[\begin{smallmatrix} 0 & -2 \\ 1 & \phantom{-}0 \end{smallmatrix}\right] \in G$.  Then 
$$
g \cdot \psi = -\left(\frac{1}{2}x^2 + 2y^2 + 2xy + 4x + 8y + 8\right)\d_x + \left(\frac{1}{4} x^2 + y^2 - xy - x + 2y + 1\right)\d_y,
$$
which we depict in Figure \ref{subfig:gpsi}.  We will evaluate the quadratic $G$-invariant functions at both $\psi$ and $g \cdot \psi$.  The only graphs in Table \ref{table:GL2} that correspond to nonvanishing invariants on $\psi$ are the finite dynamical systems, presented in rows 1, 4, and 5.  We compute these invariants in Table \ref{table:vectorfields}.
\end{example}

\begin{table}[t]
\centering
\begin{tblr}{vlines, row{3-5}={20pt}, columns={c}} \hline
Finite dynamical system $\Gamma$; &  & \\
let $I = \varphi \circ s(\Gamma)$ & $I(\psi)$ & $I(g \cdot \psi)$\\\hline[2pt]
$\v{
    \begin{tikzpicture}[>=stealth, scale=0.6]
\GraphInit[vstyle=Simple]
\tikzset{VertexStyle/.append style={minimum size=6pt, inner sep=1pt}}
\Vertex{A} \EA(A){B}
\Loop[dist=15pt, style={thick,->}](A)
\Loop[dist=15pt,dir=EA, style={thick,->}](B)
\end{tikzpicture}
}$ & 4 & 4 \\ \hline
$\v{
    \begin{tikzpicture}[>=stealth, scale=0.6]
\GraphInit[vstyle=Simple]
\tikzset{VertexStyle/.append style={minimum size=6pt, inner sep=1pt}}
\Vertex{A} \EA(A){B}
\Edge[style={thick,->, bend left}](A)(B)
\Edge[style={thick,->, bend left}](B)(A)
\end{tikzpicture}
}$ & 36 & 36\\ \hline
$\v{
    \begin{tikzpicture}[>=stealth, scale=0.6]
\GraphInit[vstyle=Simple]
\tikzset{VertexStyle/.append style={minimum size=6pt, inner sep=1pt}}
\Vertex{A} \EA(A){B}
\Edge[style={thick,->}](B)(A)
\Loop[dist=15pt, style={thick,->}](A)
\end{tikzpicture}
}$\hspace{8pt} & 16 & 16 \\\hline
   \end{tblr}
   \caption{Nonvanishing quadratic invariants on the vector fields $\psi$ and $g \cdot \psi$ from Example~\ref{example:vectorfields}.}
   \label{table:vectorfields}
\end{table}

\subsection{The orthogonal group} 

\begin{example}
We consider the familiar example of conic sections in $\mathbb R^2$.  The general form of a conic section is 
$$Ax^2 + Bxy + Cy^2 + Dx + Ey + F = 0,$$
where $A,\ldots,F \in \mathbb R$.  (We remark that $\O_n$ is defined over $\mathbb R$, and its real points constitute the subgroup $\O(n,\mathbb R)$.)  In Table \ref{table:conics}, we list the well-known invariants on conic sections, under the action of $\O(2,\mathbb R)$; the right-hand column presents the corresponding graphs from Algorithm \ref{alg}, where $G = \O_2$ is the complexification of $\O(2,\mathbb R)$.

\begin{table}[t]
\centering
\begin{tblr}{|c|c|} \hline
Classical invariant & Corresponding graph \\ \hline[2pt]
   $F$  & $\v{
    \begin{tikzpicture}[>=stealth, scale=0.4]
\GraphInit[vstyle=Simple]
\tikzset{VertexStyle/.append style={minimum size=4pt, inner sep=1pt}}
\Vertex{A}
\end{tikzpicture}
}$  \\ \hline
  $A+C$ & $\frac{1}{2}\!\v{
    \begin{tikzpicture}[>=stealth, scale=0.4]
\GraphInit[vstyle=Simple]
\tikzset{VertexStyle/.append style={minimum size=4pt, inner sep=1pt}}
\Vertex{A}
\Loop[dist=30pt, style={thick,-}](A)
\end{tikzpicture}
}$ \\ \hline
  $D^2 + E^2$ & $\v{
    \begin{tikzpicture}[>=stealth, scale=0.4]
\GraphInit[vstyle=Simple]
\tikzset{VertexStyle/.append style={minimum size=4pt, inner sep=1pt}}
\Vertex{A} \EA(A){B}
\Edge[style={thick,-}](A)(B)
\end{tikzpicture}
}$ \\ \hline
   $B^2 - 4AC$ & $\frac{1}{2}[ \v{
    \begin{tikzpicture}[>=stealth, scale=0.4]
\GraphInit[vstyle=Simple]
\tikzset{VertexStyle/.append style={minimum size=4pt, inner sep=1pt}}
\Vertex{A} \EA(A){B}
\Edge[style={thick,-, bend left}](A)(B)
\Edge[style={thick,-, bend left}](B)(A)
\end{tikzpicture}
} - \v{
    \begin{tikzpicture}[>=stealth, scale=0.4]
\GraphInit[vstyle=Simple]
\tikzset{VertexStyle/.append style={minimum size=4pt, inner sep=1pt}}
\Vertex{A} \EA(A){B}
\Loop[dist=30pt, style={thick,-}](A)
\Loop[dir=EA, dist=30pt, style={thick,-}](B)
\end{tikzpicture}
}]$ \\ \hline
   $I\coloneqq\left| \begin{smallmatrix} A & B/2 & D/2 \\ B/2 & C & E/2 \\ D/2 & E/2 & F \end{smallmatrix}\right|$ & $\frac{1}{8}\left[ \v{
    \begin{tikzpicture}[>=stealth, scale=0.4]
\GraphInit[vstyle=Simple]
\tikzset{VertexStyle/.append style={minimum size=4pt, inner sep=1pt}}
\Vertex{A} \EA(A){B} \EA(B){C}
\Edge[style={thick,-}](A)(B)
\Edge[style={thick,-}](B)(C)
\end{tikzpicture}
} + \v{
    \begin{tikzpicture}[>=stealth, scale=0.3, rotate=90]
\GraphInit[vstyle=Simple]
\tikzset{VertexStyle/.append style={minimum size=4pt, inner sep=1pt}}
\Vertices{circle}{A,B,C}
\Loop[dir = EA, dist=30pt,style={thick,-}](C)
\Loop[dir = EA, dist=30pt,style={thick,-}](B)
\end{tikzpicture}
} - \v{
    \begin{tikzpicture}[>=stealth, scale=0.3, rotate=90]
\GraphInit[vstyle=Simple]
\tikzset{VertexStyle/.append style={minimum size=4pt, inner sep=1pt}}
\Vertices{circle}{A,B,C}
\Loop[dir = WE, dist=30pt,style={thick,-}](A)
\Edge[style={thick,-}](B)(C)
\end{tikzpicture}
} - \v{
    \begin{tikzpicture}[>=stealth, scale=0.3, rotate=90]
\GraphInit[vstyle=Simple]
\tikzset{VertexStyle/.append style={minimum size=4pt, inner sep=1pt}}
\Vertices{circle}{A,B,C}
\Edge[style={thick,-,, bend left}](B)(C)
\Edge[style={thick,-,bend left}](C)(B)
\end{tikzpicture}
}\right]$ \\ \hline
   $(\text{eccentricity})^2$ & $\dfrac{2\:\sqrt{2\v{
    \begin{tikzpicture}[>=stealth, scale=0.4]
\GraphInit[vstyle=Simple]
\tikzset{VertexStyle/.append style={minimum size=4pt, inner sep=1pt}}
\Vertex{A} \EA(A){B}
\Edge[style={thick,-, bend left}](A)(B)
\Edge[style={thick,-, bend left}](B)(A)
\end{tikzpicture}
}-\v{
    \begin{tikzpicture}[>=stealth, scale=0.4]
\GraphInit[vstyle=Simple]
\tikzset{VertexStyle/.append style={minimum size=4pt, inner sep=1pt}}
\Vertex{A} \EA(A){B}
\Loop[dist=15pt, style={thick,-}](A)
\Loop[dir=EA, dist=15pt, style={thick,-}](B)
\end{tikzpicture}
}}}{\operatorname{sgn} I \!\v{
    \begin{tikzpicture}[>=stealth, scale=0.4]
\GraphInit[vstyle=Simple]
\tikzset{VertexStyle/.append style={minimum size=4pt, inner sep=1pt}}
\Vertex{A}
\Loop[dist=30pt, style={thick,-}](A)
\end{tikzpicture}
} + \sqrt{2\v{
    \begin{tikzpicture}[>=stealth, scale=0.4]
\GraphInit[vstyle=Simple]
\tikzset{VertexStyle/.append style={minimum size=4pt, inner sep=1pt}}
\Vertex{A} \EA(A){B}
\Edge[style={thick,-, bend left}](A)(B)
\Edge[style={thick,-, bend left}](B)(A)
\end{tikzpicture}
}-\v{
    \begin{tikzpicture}[>=stealth, scale=0.4]
\GraphInit[vstyle=Simple]
\tikzset{VertexStyle/.append style={minimum size=4pt, inner sep=1pt}}
\Vertex{A} \EA(A){B}
\Loop[dist=15pt, style={thick,-}](A)
\Loop[dir=EA, dist=15pt, style={thick,-}](B)
\end{tikzpicture}
}}}$ \\ \hline
   \end{tblr}
   \caption{Invariants of conic sections under the action of $\O(2,\mathbb R)$.}
\label{table:conics}
   \end{table}

\end{example}

\subsection{The symplectic group} Let $n=1$.  In this case, since $G = \Sp_2 \cong \SL_2$, our Algorithm \ref{alg} actually finds $\SL_2$-invariants on $\Psi \cong \C[x,y]$.  
Whereas early classical invariant theory focused on $\SL_2$-invariants on binary $m$-forms, i.e., on the space $S^m(\C^2)$ spanned by $\{x^m, x^{m-1}y, \ldots, xy^{m-1}, y^m\}$, our approach comprises the $\SL_2$-invariants on the direct sum of all these spaces, namely the entire polynomial ring $\C[x,y] \cong \bigoplus_{m=0}^\infty S^m(\C^2)$.

In \cite{OS}, and subsequently in \cite[Ch.~7]{olver}, Olver describes a graphical method for invariants (and more generally, covariants) of binary forms, which can be regarded as a special case of our Algorithm \ref{alg}.  Specifically, when $G = \Sp_2$ and $\Gamma$ is $m$-regular (i.e., every vertex has degree $m$), our corresponding invariant $\varphi \circ s(\Gamma)$ restricts to an invariant on the space of binary $m$-forms, which coincides (up to a constant  multiple) with the invariant obtained from the same graph in Olver's method.  (Although Olver's graphs are directed, reversing an arrow just multiplies the corresponding covariant by $-1$.
) 

\begin{remark}
In the $n=1$ case, our quadratic $r_{ij}=x_{1i}x_{2j}-x_{2i}x_{1j}$ is just the determinant of the minor corresponding to columns $i$ and $j$ of a $2 \times d$ matrix, i.e., a Pl\"ucker coordinate.  In classical invariant theory, this determinant is often denoted by $[ij]$, or more commonly by bracketed Greek letters, and historically has been named a \emph{symbolic determinantal factor} \cite{GY}, a \emph{bracket factor} \cite{weyl, olver},  or a \emph{homogenized root} \cite{KR}.
\end{remark}

\begin{example}[Invariants of binary forms]
In Table \ref{table:binaryforms}, we list the well-known fundamental invariants of the binary quadric, cubic, and quartic ($m=2,3,4$), along with their corresponding graphs from Algorithm \ref{alg}, up to a scale factor.  In place of our $\c$ notation, we write the coefficients as capital letters $A,B,\ldots$ in descending order of degree in $x$.  The symbol $\Delta_m$ denotes the discriminant of the binary $m$-form.  In the case of the quartic, we follow \cite[p.~29]{olver} in denoting the two fundamental invariants\footnote{In \cite[p.~205]{GY}, these invariants are denoted by $I$ and $J$, while their multiples are written as $i\coloneqq2I$ and $j\coloneqq6J$.} by $i$ and $j$; then the discriminant is $\Delta_4 = i^3 - 27j^2$.  Graphically, the product of invariants corresponds to the union of their corresponding graphs.  Note that the graphs in Table \ref{table:binaryforms} are indeed $m$-regular, with the number of vertices giving the degree of the invariant in the coefficients.

\begin{table}[t]
    \centering
\centerline{\begin{tblr}{colspec={|Q[c,m]|Q[c,m]|Q[c,m]|Q[c,m]|}}
\hline
Binary form & Classical invariant & Graph & Scale factor\\ \hline[2pt]
Quadric $(m=2)$ & $\Delta_2 = B^2-4AC$ &  $\v{
    \begin{tikzpicture}[>=stealth, scale=0.8]
\GraphInit[vstyle=Simple]
\tikzset{VertexStyle/.append style={minimum size=6pt, inner sep=1pt}}
\Vertex{A} \EA(A){B}
\Edge[style={thick,-, bend left}](A)(B)
\Edge[style={thick,-, bend left}](B)(A)
\end{tikzpicture}
}$ & $-2$\\ \hline
Cubic $(m=3)$ & $\Delta_3 = B^2 C^2 - 4 A C^3 - 4 B^3 D + 18 A B C D - 27 A^2 D^2$ & $\v{
    \begin{tikzpicture}[>=stealth, scale=0.8]
\GraphInit[vstyle=Simple]
\tikzset{VertexStyle/.append style={minimum size=6pt, inner sep=1pt}}
\Vertex{A} \EA(A){B} \SO(A){C} \SO(B){D}
\Edge[style={thick,-, bend left=15}](A)(B)
\Edge[style={thick,-, bend left=15}](B)(A)
\Edge[style={thick,-, bend left=15}](C)(D)
\Edge[style={thick,-, bend left=15}](D)(C)
\Edge[style={thick,-}](A)(C)
\Edge[style={thick,-}](B)(D)
\end{tikzpicture}
}$ & $96$\\ \hline
 \SetCell[r=2]{m} {Quartic $(m=4)$} & $i=AE - \frac{1}{4} BD + \frac{1}{12}C^2$ & $\v{
    \begin{tikzpicture}[>=stealth, scale=0.8]
\GraphInit[vstyle=Simple]
\tikzset{VertexStyle/.append style={minimum size=6pt, inner sep=1pt}}
\Vertex{A} \EA(A){B}
\Edge[style={thick,-, bend left}](A)(B)
\Edge[style={thick,-, bend left}](B)(A)
\Edge[style={thick,-, bend left = 10}](A)(B)
\Edge[style={thick,-, bend left = 10}](B)(A)
\end{tikzpicture}
}$ & $1152$\\ \hline
 & $j=\begin{vmatrix} A & B/4 & C/6 \\ B/4 & C/6 & D/4 \\ C/6 & D/4 & E \end{vmatrix}$ & $\v{
    \begin{tikzpicture}[>=stealth, scale=0.5, rotate=90]
\GraphInit[vstyle=Simple]
\tikzset{VertexStyle/.append style={minimum size=6pt, inner sep=1pt}}
\Vertices{circle}{A,B,C} \Edge[style={thick,-, bend left=15}](B)(A)
\Edge[style={thick,-, bend left=15}](A)(B)
\Edge[style={thick,-, bend left=15}](B)(C)
\Edge[style={thick,-, bend left=15}](C)(B)
\Edge[style={thick,-, bend left=15}](A)(C)
\Edge[style={thick,-, bend left=15}](C)(A)
\end{tikzpicture}
}$ & $82944$\\ \hline
\end{tblr}}
    \caption{Fundamental invariants of the binary quadric, cubic, and quartic.}
    \label{table:binaryforms}
\end{table}

\end{example}

\begin{example}
In Olver's method, for fixed $m$, a graph that is not $m$-regular corresponds to a \textit{covariant} on binary $m$-forms (i.e., an $\SL_2$-invariant polynomial in not just the coefficients of the form, but also in the variables $x$ and $y$).  For us, however, \textit{all} graphs correspond to true invariants on $\C[x,y]$.  Consequently, certain syzygies among Olver's graphs do not carry over into our graphs, specifically Olver's ``Rule \#2''; see illustration (7.3) in \cite[p.~135]{olver}.  As a specific example, take the following graph $\Gamma \in \mathcal G^3_4$, which for Olver in \cite[p.~135]{olver} corresponds to the zero covariant on a binary form of degree $m \geq 3$:
 $$\Gamma =  \v{
    \begin{tikzpicture}[>=stealth, rotate=-150, scale=0.6]
\GraphInit[vstyle=Simple]
\tikzset{VertexStyle/.append style={minimum size=6pt, inner sep=1pt}}
\Vertices{circle}{A,B,C}
\tikzset{EdgeStyle/.style = {-, very thick}}
\Edge(A)(C) \Edge(C)(B)
\tikzset{EdgeStyle/.style = {-, very thick, bend left}}
\Edge(A)(B) \Edge(B)(A)
\end{tikzpicture}
}$$
 For us, however, $\Gamma$ corresponds to a nontrivial degree-$3$ invariant on $\C[x,y]$:
 $$
 \varphi \circ s(\Gamma) = 16\c_{12}^2 \c^{}_{20}-8\c^{}_{11}\c^{}_{12}\c^{}_{21}-48\c^{}_{03}\c^{}_{20}\c^{}_{21}+16\c^{}_{02}\c^2_{21}+72\c^{}_{03}\c^{}_{11}\c^{}_{30}-48\c^{}_{02}\c^{}_{12}\c^{}_{30,}
 $$
 where we have suppressed the parentheses around the ordered pairs in the subscripts.  Observe that the vertices of $\Gamma$ have degrees $2$, $3$, and $3$; this is reflected in the fact that each term in the corresponding invariant takes the form $\c_\al \c_\be \c_\gamma$, where $\al,\be,\gamma \in \N^2$ and
 $$
 |\al|=2,\quad |\be|=3, \quad |\gamma|=3.
 $$
\end{example}

\section*{Appendix: Mathematica code}

\subsection*{General linear group}

\bigskip
\tiny
\begin{verbatim}
(* DEFINITIONS *)

(* Below, "chat" stands for "c_hat." *)

chat[alpha_, beta_] := Times@@Factorial@alpha * Times@@Factorial@beta * c[alpha, beta];

r[i_, j_] := Sum[y[l, i] x[l, j], {l, n}];

(* The function "s" sends a d-by-d adjacency matrix A to a product of the r[i,j]. *)

s[A_] := Product[r[i, j]^A[[i, j]], {i, d}, {j, d}] // Expand;

phi[term_] := (term /. {x[i_, j_] -> 1, y[i_, j_] -> 1}) *
   Product[chat[Table[Exponent[term, x[i, j]], {i, n}], 
                Table[Exponent[term, y[i, j]], {i, n}]], {j, d}];

(* The function "invar" sends an adjacency graph A to its corresponding invariant.  
When n=1, then s(A) has only one term. *)

invar[A_] := If[n == 1 || Total[A, 2] == 0, phi@s@A, phi/@s@A];
\end{verbatim}

\begin{verbatim}
(*BEGIN ALGORITHM HERE.*)

(* Input: n-value.
Input: the d-by-d adjacency matrix A of a directed graph.
Output: the graph and its associated invariant .*)

n = 2;
A = {{1,0}, {0,1}};
d = Length@A; i = invar@A;
{AdjacencyGraph[A, DirectedEdges -> True], 
i /. c[list1_, list2_] :> Subscript[c, {list1, list2}]}
\end{verbatim}

\subsection*{Orthogonal group} 

\begin{verbatim}
(* DEFINITIONS *)

chat[alpha_] := Times@@Factorial@alpha*c[alpha];

r[i_, j_] := Sum[x[l, i] x[l, j], {l, n}];

s[A_] := Product[r[i, j]^A[[i, j]], {i, d}, {j, i, d}] // Expand;

phi[term_] := (term /. x[i_, j_] -> 1) *
   Product[chat@Table[Exponent[term, x[i, j]], {i, n}], {j, d}];

invar[A_] := 
  If[n == 1 || Total[PadLeft@A, 2] == 0, phi@s@PadLeft@A, phi/@s@PadLeft@A];

(* "TriToSym" returns a d-by-d symmetric matrix, 
given the upper-triangular entries "upper"
as a list of lists of lengths d, d-1, ... , 1. *)

TriToSym[upper_] := 
 PadLeft[upper] + Transpose[UpperTriangularize[PadLeft[upper], 1]];

(* BEGIN ALGORITHM HERE *)

(* Input: n-value.
Input: the upper-triangular entries U of a symmetric d-by-d adjacency matrix.
Output: the graph and its associated invariant .*)

n = 2; 
U = {{1, 2, 3}, {4, 5}, {6}};
d = Length@U; i = invar@U;
{AdjacencyGraph@TriToSym@U, i /. c[list_] :> Subscript[c, list]}
\end{verbatim}

\subsection*{Symplectic group}

\begin{verbatim}
(* DEFINITIONS *)

chat[alpha_] := Times@@Factorial@alpha*c[alpha];

r[i_, j_] := Sum[x[l + n, i] x[l, j] - x[l, i] x[l + n, j], {l, n}];

s[A_] := Product[r[i, j]^A[[i, j]], {i, Length@A - 1}, {j, i + 1, Length@A}] // Expand;

phi[term_] := (term /. x[i_, j_] -> 1) *
   Product[chat@Table[Exponent[term, x[i, j]], {i, 2 n}], {j, d}];

(* "StrictTriToSkew" returns a d-by-d skew-symmetric matrix,
given the strictly upper-triangular entries "upper"
as a list of lists of lengths d-1, ..., 1. *)

StrictTriToTri[upper_] := PadLeft[Append[Prepend[0] /@ upper, {0}]];
StrictTriToSkew[upper_] := StrictTriToTri@upper - Transpose@StrictTriToTri@upper;

invar[upper_] := 
  If[Total[upper, 2] == 0, phi@s@StrictTriToSkew@upper, phi/@s@StrictTriToSkew@upper];

(*BEGIN ALGORITHM HERE *)

(* Input: n-value.
Input: the strictly upper-triangular entries U of a d-by-d adjacency matrix.
Output: the graph and its associated invariant. *)

n = 1;
U = {{2}};
d = Length@U + 1; i = invar@U;
{AdjacencyGraph[UpperTriangularize@StrictTriToSkew@U + 
   Transpose[UpperTriangularize@StrictTriToSkew@U]],
   i /. c[list_] :> Subscript[c, list]}
\end{verbatim}

\normalsize

\subsection*{Acknowledgments}

We would like to thank the referee for the careful reading and for the many helpful suggestions that improved the quality of this paper.
We are especially grateful to the referee for pointing out the connection to the (largely forgotten) classical notion of perpetuants, which in turn led us to find the recent paper~\cite{KP}.

\begin{Backmatter}

\bibliographystyle{alpha}

\bibliography{references}

\begin{thebibliography}{HSW21}

\bibitem[Bur74]{Burge}
William~H. Burge.
\newblock Four correspondences between graphs and generalized {Y}oung tableaux.
\newblock {\em J. Combinatorial Theory Ser. A}, 17:12--30, 1974.

\bibitem[Con94]{Conca94}
Aldo Conca.
\newblock Gr\"{o}bner bases of ideals of minors of a symmetric matrix.
\newblock {\em J. Algebra}, 166(2):406--421, 1994.

\bibitem[Cvi08]{Cvitanovic}
Predrag Cvitanovi\'{c}.
\newblock {\em Group theory: birdtracks, {L}ie's, and exceptional groups}.
\newblock Princeton University Press, Princeton, NJ, 2008.

\bibitem[Dol03]{dolgachev}
Igor Dolgachev.
\newblock {\em Lectures on invariant theory}, volume 296 of {\em London
  Mathematical Society Lecture Note Series}.
\newblock Cambridge University Press, Cambridge, 2003.

\bibitem[DRS74]{DRS}
Peter Doubilet, Gian-Carlo Rota, and Joel Stein.
\newblock On the foundations of combinatorial theory. {IX}. {C}ombinatorial
  methods in invariant theory.
\newblock {\em Studies in Appl. Math.}, 53:185--216, 1974.

\bibitem[GW09]{gw}
Roe Goodman and Nolan~R. Wallach.
\newblock {\em Symmetry, representations, and invariants}, volume 255 of {\em
  Graduate Texts in Mathematics}.
\newblock Springer, Dordrecht, 2009.

\bibitem[GY10]{GY}
John~Hilton Grace and Alfred Young.
\newblock {\em The algebra of invariants}.
\newblock Cambridge Library Collection. Cambridge University Press, Cambridge,
  2010.
\newblock Reprint of the 1903 original.

\bibitem[HSW21]{hsw}
Alexander Heaton, Songpon Sriwongsa, and Jeb~F. Willenbring.
\newblock Branching from the general linear group to the symmetric group and
  the principal embedding.
\newblock {\em Algebr. Comb.}, 4(2):189--200, 2021.

\bibitem[Kem86]{kempe}
Alfred~B. Kempe.
\newblock On the application of {C}lifford's graphs to ordinary binary
  quantics.
\newblock {\em Proc. Lond. Math. Soc.}, 17:107--121, 1885/86.

\bibitem[Knu70]{Knuth}
Donald~E. Knuth.
\newblock Permutations, matrices, and generalized {Y}oung tableaux.
\newblock {\em Pacific J. Math.}, 34:709--727, 1970.

\bibitem[KP21]{KP}
Hanspeter Kraft and Claudio Procesi.
\newblock Perpetuants: a lost treasure.
\newblock {\em Int. Math. Res. Not. IMRN}, 5:3597--3632, 2021.

\bibitem[KR84]{KR}
Joseph P.~S. Kung and Gian-Carlo Rota.
\newblock The invariant theory of binary forms.
\newblock {\em Bull. Amer. Math. Soc. (N.S.)}, 10(1):27--85, 1984.

\bibitem[Kup96]{Kuperberg}
Greg Kuperberg.
\newblock Spiders for rank-2 {L}ie algebras.
\newblock {\em Comm. Math. Phys.}, 180(1):109--151, 1996.

\bibitem[Olv99]{olver}
Peter~J. Olver.
\newblock {\em Classical invariant theory}, volume~44 of {\em London
  Mathematical Society Student Texts}.
\newblock Cambridge University Press, Cambridge, 1999.

\bibitem[OS89]{OS}
Peter~J. Olver and Chehrzad Shakiban.
\newblock Graph theory and classical invariant theory.
\newblock {\em Adv. Math.}, 75(2):212--245, 1989.

\bibitem[Pen05]{Penrose}
Roger Penrose.
\newblock {\em The road to reality: a complete guide to the laws of the
  universe}.
\newblock Alfred A. Knopf, Inc., New York, 2005.

\bibitem[Pro07]{procesi}
Claudio Procesi.
\newblock {\em Lie groups: an approach through invariants and representations}.
\newblock Universitext. Springer, New York, 2007.

\bibitem[Sha94]{PV}
Igor~R. Shafarevich, editor.
\newblock {\em Algebraic geometry. {IV}}, volume~55 of {\em Encyclopaedia of
  Mathematical Sciences}.
\newblock Springer-Verlag, Berlin, 1994.

\bibitem[Sta99]{stanley2}
Richard~P. Stanley.
\newblock {\em Enumerative combinatorics. {V}ol. 2}, volume~62 of {\em
  Cambridge Studies in Advanced Mathematics}.
\newblock Cambridge University Press, Cambridge, 1999.

\bibitem[Stu90]{SturmfelsGB}
Bernd Sturmfels.
\newblock Gr\"{o}bner bases and {S}tanley decompositions of determinantal
  rings.
\newblock {\em Math. Z.}, 205(1):137--144, 1990.

\bibitem[Stu08]{sturmfels}
Bernd Sturmfels.
\newblock {\em Algorithms in invariant theory}.
\newblock Texts and Monographs in Symbolic Computation. SpringerWienNewYork,
  Vienna, second edition, 2008.

\bibitem[Syl78]{Sylvester}
James~Joseph Sylvester.
\newblock Chemistry and algebra.
\newblock {\em Nature}, 17(284), 1878.

\bibitem[Syl82]{SylvesterSubvariants}
James~Joseph Sylvester.
\newblock On subvariants, i.e. semi-invariants to binary quantics of an
  unlimited order.
\newblock {\em Amer. J. Math.}, 5(1-4):79--136, 1882.

\bibitem[Wal17]{WallachGIT}
Nolan~R. Wallach.
\newblock {\em Geometric invariant theory over the real and complex numbers}.
\newblock Universitext. Springer, Cham, 2017.

\bibitem[Wey39]{weyl}
Hermann Weyl.
\newblock {\em The classical groups: their invariants and representations}.
\newblock Princeton University Press, Princeton, N.J., 1939.

\end{thebibliography}

\printaddress

\end{Backmatter}

\end{document}